\documentclass[12 pt]{article}
\usepackage[utf8]{inputenc}
\usepackage{amsfonts}
\usepackage{amsthm}
\usepackage{tikz}
\usepackage{pgfplots}
\pgfplotsset{compat=1.12}
\usepackage{treeFile}
\usepackage{graphicx}
\usepackage{mathtools}
\usepackage{extarrows}
\usepackage{hyperref}
\usepackage{listings}
\usepackage{pythonhighlight}

\newtheorem*{remark}{Remark}

\newcommand{\bR}{\mathbb{R}}
\newcommand{\bC}{\mathbb{C}}

\newcommand{\cB}{\mathcal{B}}
\newcommand{\cC}{\mathcal{C}}
\newcommand{\cD}{\mathcal{D}}
\newcommand{\cP}{\mathcal{P}}
\newcommand{\cJ}{\mathcal{J}}
\newcommand{\cM}{\mathcal{M}}
\newtheorem{definition}{Definition}
\newtheorem{example}{Example}

\title{Deformation Spaces and Static Animations}
\author{Gabriel Dorfsman-Hopkins}
\date{\today}

\begin{document}
\maketitle

\begin{abstract}
      We study applications of 3D printing to the broad goal of understanding how mathematical objects vary continuously in families.  To do so, we model the varying parameter as the vertical axis of a 3D print, introducing the notion of a \textit{static animation}: a 3D printed object each of whose layers is a member of the continuously deforming family.  We survey examples and draw connections to algebraic geometry, complex dynamics, chaos theory, and more.  We also include a detailed tutorial (with accompanying code and files) so that the reader can create static animations of their own.
\end{abstract}

\section{Introduction}
Across mathematical disciplines there is considerable interest in understanding how a mathematical object can deform according to a shift in a single parameter, producing a varying family of related mathematical objects.  If these objects have small enough dimension, one can vary the parameter over time and display how the object changes as an animation on a screen. With a 3D printer we can replace the time dimension with the vertical $z$-axis, and stack each successively deformed layer on top of the last. This turns the entire parametrized family of objects into a single solid which we can hold in our hands, giving an entire new physical dimension to the deformation.\\

The study of continuously varying families is ubiquitous in mathematics.  To name a few examples:  In algebraic geometry, one studies degenerations of curves, surfaces, and higher dimensional algebraic varieties, putting them in continuously varying \textit{flat families} which interpolate between one space and the next \cite[Chapter III.9]{Hartshorne}.  In complex dynamics, one asks about limiting behavior when iterating a holomorphic function $f$, and compares this to the limiting behavior of a slight perturbation of $f$, to dazzling effects \cite{Milnor}.  In topology, the question of whether one space can continuously deform to another forms the notion of a homotopy, defining the fundamental invariant in the rich and expansive field of homotopy theory \cite[Chapter 4]{Hatcher}.  In the study of chaotic dynamical systems, one explores how slight changes in initial conditions can radically affect longer term outcomes \cite{strogatz}.  In fact, one would be hard pressed to find a field of mathematics which does not contain a multitude of problems that aim to understand how a mathematical object or system varies when perturbed.\\

In this article we propose a general framework which can be used to explore situations like those described in the previous paragraph.  Namely, we describe ways to unify an entire family of continuously varying objects into a single mathematical object or \textit{deformation space}.  When the members of this family are 1 or 2 dimensional, we argue that 3D printing this deformation space provides new perspectives on the illustration of continuously deforming geometry, beyond that which an animation on a screen can achieve.  We call these 3D models \textit{static animations}.\\

There are two main components to this article, one more theoretical and one more interactive.  The theoretical component precisely formulates the notion of a static animation and describes connections to various mathematical areas, showcasing examples and suggesting avenues for exploration.  The interactive component provides detailed instructions for the reader to create static animations of their own, and is accompanied by code and files in a public github repository \cite{github}.\\

The structure of the article is therefore as follows.  In Section \ref{Theory} we give a rigorous definition of the notion of a static animation, and illustrate the concept with an example of a varying family of polar flowers.  In Section \ref{treeSection} we introduce the author's first 3D printed deformation space, which illustrates a family of Pythagorean tree fractals.  In Section \ref{juliaSection} we enter the world of complex dynamics, showcasing various prints of deformations of Julia sets, including some our own work, as well as work Bernat Espigul\'e, Caroline Davis, and Bethany Mussman.\\

We then arrive at the interactive portion in Section \ref{How?}, describing a procedure to create static animations that was developed with Bernat Espigul\'e.  This workflow can be understood as a sort of \textit{reverse tomography}.  Tomography is a way to visualize the internal structure of 3D object by producing a collection 2D pictures of cross sections at various depths.  For example, an MRI is a type of tomography which can be applied to the human body.  The process we describe goes in the other direction, starting with a collection of 2D images and stitching them together into a 3D object.\\

The upshot of this particular workflow is that the user can create static animations using mostly 2D modelling methods$-$minimizing the need for more difficult 3D modeling software.  The main idea is to feed a stack of 2D images to an open source medical and molecular visualization software called Chimera \cite{chimera},\cite{chimeraSoftware} which will stitch them together into a 3D model.\footnote{As far as the author is aware, this functionality was initially developed to reconstruct 3D models from medical tomographical images such as MRI images.}  After broadly describing the workflow, we also include a detailed step-by-step tutorial allowing the reader to follow along and create a static animation of their own.  Accompanying the tutorial we have our own code and files from each step of the process, including a file of the finalized 3D model.  It is the hope of the author that after working through this tutorial, any reader able to make 2D animations on a screen will also be able to make 3D static animations.\\

The idea of viewing a continuously deforming family of a geometric objects as a geometric object in its own right is certainly not new, and in fact it was the algebrogeometric perspective of viewing so-called flat morphisms as parameter spaces which motivated the author's exploration into this field.  We provide more details and an example of this connection in the appendix.

\section*{Acknowledgments}
The author's initial exploration into static animations was as the fabrication manager of the Washington Experimental Math Lab at the University of Washington$-$supported by NSF grant DMS-1559860$-$and the author thanks Jayadev Athreya for the technology and support in exploring math from an experimental and illustrative perspective during this time.  The author thanks Bernat Espigul\'e for his generosity and collaboration in creating static animations of Julia sets, Caroline Davis for extremely enlightening conversations, and Edmund Harriss for taking many of the photographs of their models which are included in this manuscript.  We also thank the anonymous referee for many helpful comments and contributions, including some of the python code in our tutorial section, as well as the description of our process as a \textit{reverse tomography}.  Much of the development of the ideas and models were made during the 2019 ICERM semester titled \textit{Illustrating Mathematics}$-$supported by NSF grant DMS-1439786$-$and the author thanks all the participants for their collaboration and support during that semester.  The author also thanks the participants of the 2021 Illustrating Mathematics Summer School at PCMI, and the 2022 Joint Mathematics Meetings mini-course \textit{3D Printing: Challenges and Applications}.  Many of the 3D models were made with UCSF Chimera, developed by the Resource for Biocomputing, Visualization, and Informatics at the University of California, San Francisco, with support from NIH P41-GM103311.  The author is currently supported by NSF grant DMS-1646385 as part of the Research Training Group in Arithmetic Geometry at UC Berkeley.
%The author also thanks Daniel Piker for developing and sharing an approach towards making static animations related to the double pendulum experiment.

\section{Static Animations}\label{Theory}
In what follows we give a precise mathematical formulation of what we mean when referring to a static animation.
\begin{definition}[Static Animations]\label{maindef}
    Suppose we are given some interval $I\subseteq\mathbb{R}$, and for each $t\in I$ we are given a subset $B_t\subseteq\mathbb{R}^2$.  The \textit{static animation} associated to the data $(I,(B_t)_{t\in I})$ is the following subset of $\bR^3$.
    \[\cB = \{(b,t):t\in I,b\in B_t\}\subseteq\bR^3.\]
    The sets $B_t$ are called the \textit{frames} of the static animation.
\end{definition}
This gives us some subset of $\bR^3$, which we could hope to 3D print.\footnote{That said, if the set $B_t$ varies too wildly with $t$, we will probably be unable to do so.  To have any hope, we'd like $B_t$ to vary continuously in $t$.  Loosely speaking, this means that if $t\approx t'$ then $B_t\approx B_{t'}$.  This can be made precise in a number of ways, depending on the context.}  Observe that a static animation comes equipped with a projection map $\pi:\cB\to I$ given by the rule $\pi(b,t) = t$.  For any $t_0\in I$ we can recover the frame $B_{t_0}$ as the the fiber of $\pi$ over $t_0$:
\[\pi^{-1}(t_0) = \{(b,t_0)\in\cB\} \cong B_{t_0}.\]
We illustrate the notion of a static animation and its associated projection map with an example of a continuously varying family of curves, interpolating between the polar flowers $B_{t_1},B_{t_2},$ and $B_{t_3}$, below.
\begin{center}
\begin{tikzpicture}
  \draw[domain=0:2*pi,scale=1,samples=500] plot ({deg(\x)}:{2 + cos(5 * (\x r))});
  \node at (-3,2) {$B_{t_1}$};
  \begin{scope}[shift={(3.5,-6)}]
    \draw[domain=0:2*pi,scale=1,samples=500] plot ({deg(\x)}:{2 + .2*cos(5 * (\x r) + 100)});
    \node at (-2,2) {$B_{t_3}$};
  \end{scope}
  \begin{scope}[shift={(7,0)}]
  \draw[domain=0:2*pi,scale=1,samples=500] plot ({deg(\x)}:{2 + .6*cos(5 * (\x r ) + 50)});
  \node at (-2.5,2) {$B_{t_2}$};
  \end{scope}
\end{tikzpicture}
\end{center}
We start by aligning these curves vertically, with the polar curve $B_t$ curve living on the plane $z=t$.
\begin{center}
\begin{tikzpicture}
\begin{axis}[hide axis]
\addplot3[variable=t,domain=0:360,samples=65] ({(2 + cos(5 * (t) ))*sin(t)},{(2 + cos(5 * (t)))*cos(t)}, 0);
\addplot3[variable=t,domain=0:360,samples=65] ({(2 + .6*cos(5 * (t) + 50))*sin(t)},{(2 + .6*cos(5 * (t)+50))*cos(t)}, .4);
\addplot3[variable=t,domain=0:360,samples=65] ({(2 + .2*cos(5 * (t)+100))*sin(t)},{(2 + .2*cos(5 * (t)+100))*cos(t)}, .8);
\end{axis}
\node at (-.5,4.25) {$B_{t_3}$};
\node at (-.5,2.75) {$B_{t_2}$};
\node at (-.5,1.25) {$B_{t_1}$};
\end{tikzpicture}
\end{center}
Then we interpolate by filling in the rest of the $t$ values to get a surface, $\cB$.  We will explicitly describe this surface in Section \ref{preciseTheory} below.
\begin{center}
\begin{tikzpicture}
    \begin{axis}[hide axis]
    \addplot3[surf,samples=50,samples y =30, domain=0:360,domain y=0.2:1] ({(2 + y*cos(5*(x)+  125*(1-y)))*sin(x)},{(2 + y*cos(5*(x)+125*(1-y)))*cos(x)}, {1-y});
    \end{axis}
\node at (-.5,2.5) {$\cB$};
\end{tikzpicture}
\end{center}
The projection $\pi:\cB\to\bR$ consists of simply plucking out the $z$-coordinate of a point on this surface.  For $t\in I$, the fiber over $t$ of this projection, $\pi^{-1}(t)$, is therefore the horizontal cross section given by intersecting $\cB$ with the plane $z=t$.  By construction, this is precisely the polar flower $B_t$, floating at height $t$.\\

We can now cap this surface to make a solid which can be 3D printed as a deformation space illustrating how to continuously interpolate between polar flowers.
\begin{center}
    \includegraphics[width=.35\paperwidth]{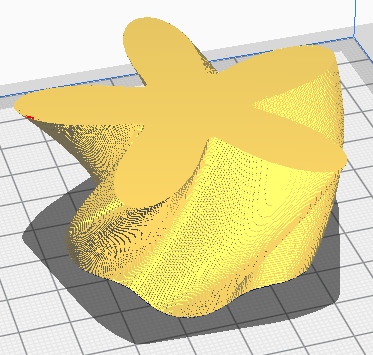}\footnote{This is a screenshot of the STL file opened and rendered in Ultimaker-Cura \cite{Cura}}\\
    \vspace{10pt}
    \includegraphics[width=.3\paperwidth]{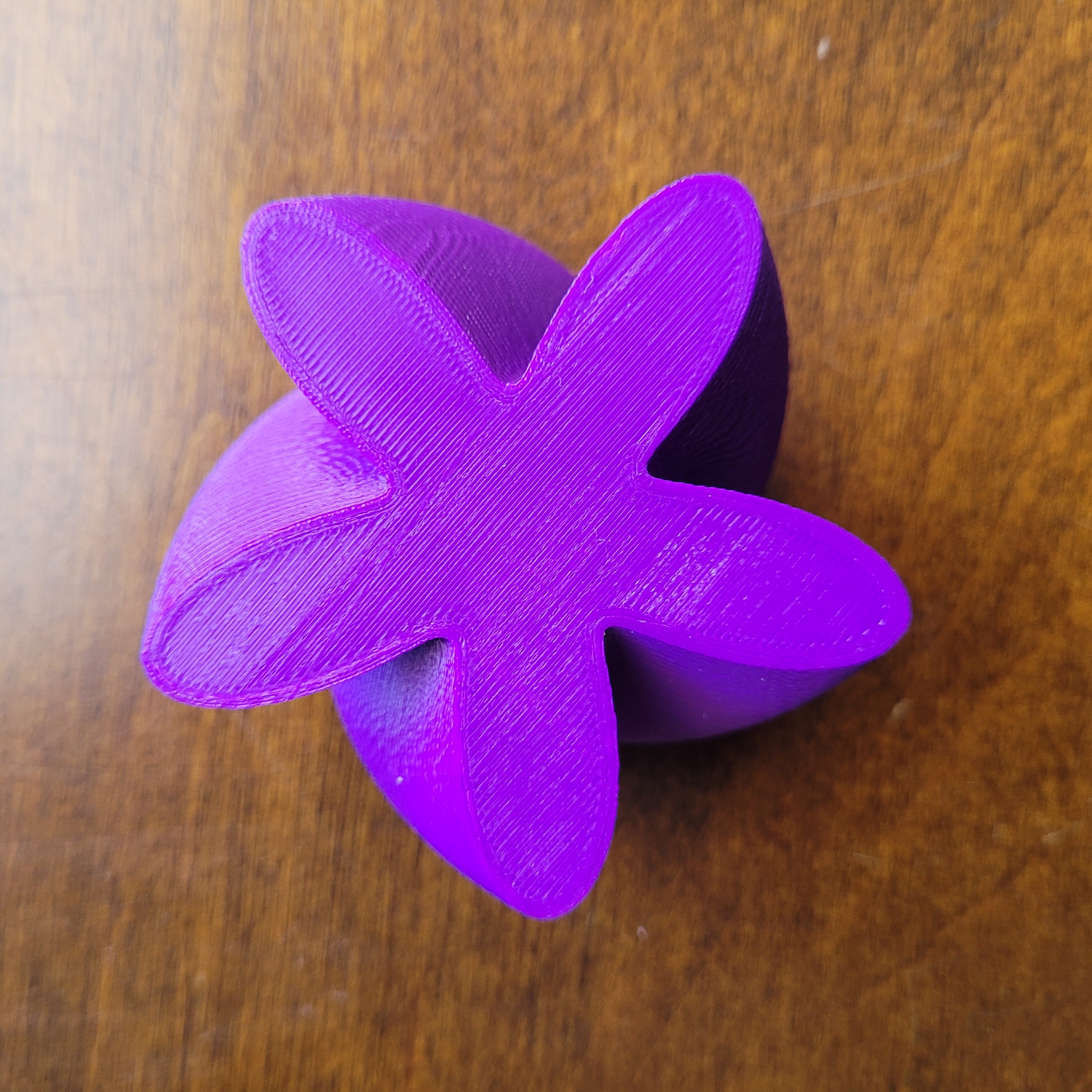}\hspace{10pt}\includegraphics[width=.3\paperwidth]{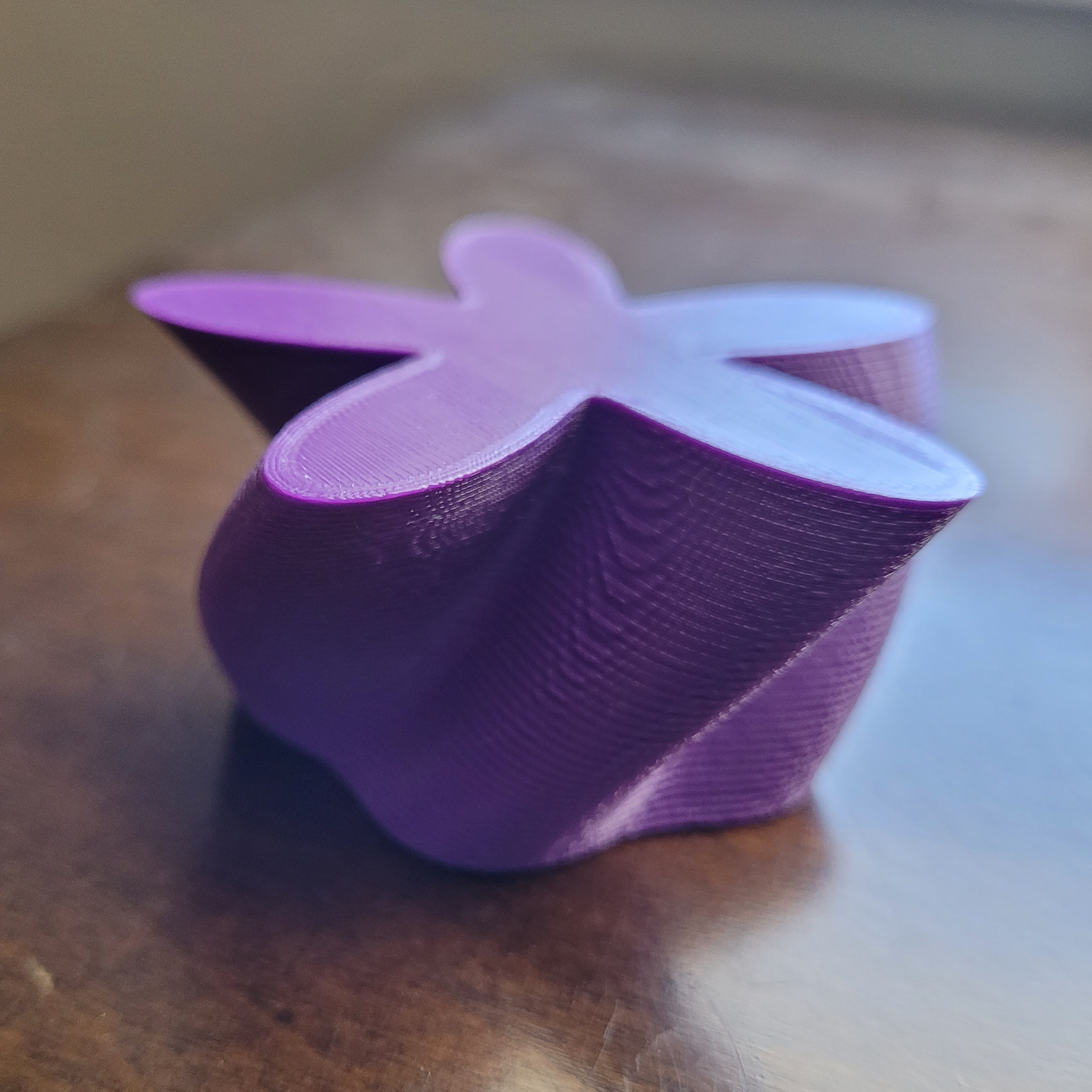}
\end{center}
We give a step by step tutorial for the creation of this object in Section \ref{tutorial}.  An STL file of the 3D model is also available at our public github repository \cite{github}, and can be downloaded to print.\\

Of course, one could also illustrate the continuous deformation of polar flowers as an animation on a screen.\footnote{We made such an animation.  Find it here: \url{http://www.gabrieldorfsmanhopkins.com/digital_fabrication/polarFlowerDeformation/flowerAnimation.gif}}  Indeed, the horizontal cross sections of the solid described above correspond precisely to the frames of this animation.  This is why we call these objects \textit{static animations}, and why we call their horizontal cross sections \textit{frames}.\\

\subsection{A Mathematical Description of The Deformation of Polar Flowers}\label{preciseTheory}
We continue with a more precise formulation of the example above.  Adopting the same notation, $B_{t_1}$ is the polar flower given by the equation
\[r = 2 + cos(5\theta)\]
in polar coordinates.  Similarly, $B_{t_3}$ is given by
\[r = 2 + \frac{1}{5}\cos(5\theta + 2\pi/5).\]
In particular, both consist of waves on the circle of radius 2, although they differ in amplitude and a shift of where the wave begins.  These correspond to two varying parameters in the equation: the amplitude is controlled by the coefficient in front of the cosine, and the shift corresponds to adding a constant term to the argument of the cosine.  These can be adjusted simultaneously to interpolate continuously  between $B_{t_1}$ and $B_{t_3}$, giving a family of curves whose general element $B_t$ is given by the equation
\[r = 2 + t\cos(5\theta + 2\pi t).\]
Indeed, letting t=1 we recover $B_{t_1}$ and letting $t=1/5$ we recover $B_{t_3}.$\footnote{We can also recover $B_{t_2}$ by letting $t=3/5$.}  From this we derive an explicit parametrization for our deformation surface
\[\rho:[0,2\pi]\times[.2,1]\longrightarrow\bR^3.\]
Giving target cylindrical coordinates, $\rho$ can be expressed by the rule
\[\rho(\theta,t) = (2+t\cos(5\theta+2\pi t),\theta,t).\]
The static animation $\cB$ is precisely the image of $\rho$, and by restricting to $\cB$ the projection $\bR^3\to\bR$ onto the $z$-axis, we recover the projection $\pi:\cB\to I$.  For $t_0\in I$ we compute that
\[\pi^{-1}(t_0) = \{\rho(\theta,t_0):\theta\in[0,2\pi]\} = \{(2+t_0\cos(5\theta+2\pi t_0),\theta,t_0):\theta\in[0,2\pi]\}.\]
This is precisely the polar flower given by $r = 2+t_0\cos(t\theta+2\pi t_0)-$that is, the frame $B_{t_0}-$floating at the fixed height $t_0$.\\

\subsection{Higher Dimensional Stacks of Frames}\label{HDF}
It is worth returning to the observation that there were two varying parameters in our family of polar flowers, one changing amplitude and another controlling the phase shift.  We controlled both using the single parameter $t$, but of course each can move independantly, say, with parameters $s$ and $t$.  This would give a 2-dimensional parameter space of polar flowers, say $\bR^2$, and for every $(s,t)\in\bR^2$ we could consider the curve
\begin{equation}\label{2DFrame}
C_{s,t}:\hspace{10pt} r = 2 + s\cos(5\theta + 2\pi t).
\end{equation}
This puts us in a situation similar to Definition \ref{maindef}, except that we have a 2-dimensional stack of frames.  Therefore, trying to build the static animation would put us in 4-space,
\[{\cC} = \{(b,s,t):b\in B_{s,t}\}\subseteq \bR^2\times\bR^2 = \bR^4.\]
${\cC}$ can also be explicitly parametrized, via the rule:
\[(\theta,s,t)\mapsto (r,\theta, z,w) = (2 + s\cos(5\theta + 2\pi t),\theta,s,t)\subseteq\bR^2\times\bR^2,\]
where the first $\bR^2$ is given polar coordinates.  Of course, this attempted static animation is now a 3-dimensional hypersurface in $\bR^4$, making it impossible to 3D print.  This is exactly the type of situation studied in Section \ref{juliaSection} where we look at deformations of Julia Sets which have a 2-dimensional parameter space (the Mandelbrot set), so rather than just give up here, we will briefly describe a way to cut down a dimension.\\

Let $I$ be an interval and consider a path $\gamma:I\to\bR^2$, which we express coordinatewise as $\gamma(t) = (x(t),y(t))$.   Then for every $t\in I$ we can consider the frame $D_t:=C_{x(t),y(t)}$ defined as in Equation (\ref{2DFrame}).  We are now back in the setup of Definition \ref{maindef} and can form the static animation
\[\cD_\gamma:=\{(d,t):t\in I, d\in D_t = C_{x(t),y(t)}\}\subseteq\bR^3.\]
For example, if we let $\gamma:[.2,1]\to\bR^2$ be the path $\gamma(t) = (1-t,1-t)$ (consisting of travelling from $(1,1)$ to $(.2,.2)$ in a straight line along the diagonal), then $\cD_\gamma$ recovers the static animation $\cB$ from Section \ref{preciseTheory}.\\

The key takeaway is the following: if one has a continuously varying family of 2D objects which vary over a higher dimensional parameter space $P$, one can still create 3D static animations by considering one parameter variations associated to paths within $P$.
\section{Deforming the Pythagorean Tree Fractal}\label{treeSection}
The Pythagorean tree fractal is a recursive plane fractal, constructed as follows.  First, one attaches the the hypotenuse of a right triangle to the side of a square.
\begin{center}
\begin{tikzpicture}[scale=.65]
\draw (0,0) -- (0,3) -- (3,3) -- (3,0) -- cycle;
\node at (4,1.5) {\Huge$\Longrightarrow$};
\draw (5,0) -- (5,3) -- (8,3) -- (8,0) -- cycle;
\draw (5,3) -- (6.5,4.5) -- (8,3) -- cycle;
\end{tikzpicture}
\end{center}
Then two squares are attached to the legs of this triangle, and the process repeats, with the hypotenuse of a right triangle affixed to the opposite side of each new square.
\begin{center}
\begin{tikzpicture}[scale=.65]
\begin{scope}[shift={(-7,0)}]
\draw (5,0) -- (5,3) -- (8,3) -- (8,0) -- cycle;
\draw (5,3) -- (6.5,4.5) -- (8,3) -- cycle;
\draw (5,3) -- (6.5,4.5) -- (5,6) -- (3.5,4.5) -- cycle;
\draw (6.5,4.5) -- (8,6) -- (9.5,4.5) -- (8,3) -- cycle;
\end{scope}
\node at (4,1.5) {\Huge$\Longrightarrow$};
\begin{scope}[shift = {(2,0)}]
\draw (5,0) -- (5,3) -- (8,3) -- (8,0) -- cycle;
\draw (5,3) -- (6.5,4.5) -- (8,3) -- cycle;
\draw (5,3) -- (6.5,4.5) -- (5,6) -- (3.5,4.5) -- cycle;
\draw (6.5,4.5) -- (8,6) -- (9.5,4.5) -- (8,3) -- cycle;
\draw (3.5,4.5) -- (3.5,6) -- (5,6);
\draw (8,6) -- (9.5,6) -- (9.5,4.5);
\end{scope}
\end{tikzpicture}
\end{center}
Squares are added to the legs of these triangles, and the process continues recursively forever.
\begin{center}
\begin{tikzpicture}
\begin{scope}[scale=1]
\PythagorasTree{3}{45}
\begin{scope}[shift={(6,0)}]
\PythagorasTree{4}{45}
\end{scope}
\end{scope}
\end{tikzpicture}
\end{center}
\begin{center}
\begin{tikzpicture}
\begin{scope}[scale = 1]
\PythagorasTree{5}{45}
\begin{scope}[shift={(6,0)}]
\PythagorasTree{6}{45}
\end{scope}
\end{scope}
\end{tikzpicture}
\end{center}
\begin{center}
\begin{tikzpicture}
\begin{scope}[scale = 1]
\PythagorasTree{7}{45}
\begin{scope}[shift={(6,0)}]
\PythagorasTree{8}{45}
\end{scope}
\end{scope}
\end{tikzpicture}
\end{center}
Notice all of the triangles in the construction above were congruent, in particular they were isosceles right triangles with acute angles of $45^\circ$.  One could ask what happens if one used, say, a 30-60-90 triangle at each step instead.  As before, one begins by aligning the hypotenuse of a 30-60-90 triangle with the top of the square.
\begin{center}
\begin{tikzpicture}[scale=.65]
\draw (0,0) -- (0,3) -- (3,3) -- (3,0) -- cycle;
\node at (4,1.5) {\Huge$\Longrightarrow$};
\draw (5,0) -- (5,3) -- (8,3) -- (8,0) -- cycle;
\draw (5,3) -- (5.75,4.29) -- (8,3) -- cycle;
\node at (6.9,3.25) {\small$30$};
\end{tikzpicture}
\end{center}
Then, following a similar process, two squares are attached to the legs of this triangle, although in this case they will not be the same size, instead matching the two differently sized legs.  Continuing, the hypotenuse of an appropriately sized 30-60-90 triangle is attached to the opposite side of each square, making sure to maintain consistent orientations.
\begin{center}
\begin{tikzpicture}[scale=.65]
\draw (0,0) -- (0,3) -- (3,3) -- (3,0) -- cycle;
\draw (0,3) -- (.75,4.29) -- (3,3) -- cycle;
\node at (1.9,3.25) {\small$30$};
\begin{scope}[shift={(0,3)},rotate=60,scale=.5]
\draw (0,0) -- (0,3) -- (3,3) -- (3,0) -- cycle;
\end{scope}
\begin{scope}[shift={(3,3)},rotate=60,scale=.866]
\draw (0,0) -- (0,3) -- (3,3) -- (3,0) -- cycle;
\end{scope}
\node at (6,1.5) {\Huge$\Longrightarrow$};
\begin{scope}[shift={(9,0)}]
\draw (0,0) -- (0,3) -- (3,3) -- (3,0) -- cycle;
\draw (0,3) -- (.75,4.29) -- (3,3) -- cycle;
\node at (1.9,3.25) {\small$30$};
\begin{scope}[shift={(0,3)},rotate=60,scale=.5]
\draw (0,0) -- (0,3) -- (3,3) -- (3,0) -- cycle;
\draw (0,3) -- (.75,4.29) -- (3,3) -- cycle;
\node at (1.9,3.3) {\tiny$30$};
\end{scope}
\begin{scope}[shift={(.75,4.29)},rotate=-30,scale=.866]
\draw (0,0) -- (0,3) -- (3,3) -- (3,0) -- cycle;
\draw (0,3) -- (.75,4.29) -- (3,3) -- cycle;
\node at (1.9,3.25) {\tiny$30$};
\end{scope}
\end{scope}
\end{tikzpicture}
\end{center}
As before, we continue this process recursively, this time observing that the global effect is somewhat different than when we used isosceles triangles.
\begin{center}
\begin{tikzpicture}
\begin{scope}[scale=1]
\PythagorasTree{3}{60}
\begin{scope}[shift={(6,0)}]
\PythagorasTree{4}{60}
\end{scope}
\end{scope}
\end{tikzpicture}
\end{center}
\begin{center}
\begin{tikzpicture}
\begin{scope}[scale = 1]
\PythagorasTree{5}{60}
\begin{scope}[shift={(6,0)}]
\PythagorasTree{6}{60}
\end{scope}
\end{scope}
\end{tikzpicture}
\end{center}
\begin{center}
\begin{tikzpicture}
\begin{scope}[scale = 1]
\PythagorasTree{7}{60}
\begin{scope}[shift={(6,0)}]
\PythagorasTree{8}{60}
\end{scope}
\end{scope}
\end{tikzpicture}
\end{center}
A right triangle is determined (up to congruence), by its smallest angle $\theta$.  In particular, for each angle $0<\theta\le45$, one can follow this construction using a right triangle whose acute angles are $\theta$ and $90-\theta$.  By varying the angle $\theta$, one obtains a varying family of planar fractals parametrized by $\theta$.  Below we display the resulting fractal for some values of $\theta$ between $45$ and $30$.  The images seem to indicate that these fractals interpolate between the two that we have seen so far, suggesting the creation of a (static) animation.
\begin{center}
    \begin{tikzpicture}
    \begin{scope}[scale=1]
    \PythagorasTree{9}{45}
    \node at (.55,-.7) {$\theta=45$};
    \begin{scope}[shift={(6,0)}]
    \PythagorasTree{9}{50}
    \node at (.55,-.7) {$\theta=40$};
    \end{scope}
    \begin{scope}[shift = {(0,-6)}]
    \PythagorasTree{9}{55}
    \node at (.55,-.7) {$\theta=35$};
    \end{scope}
    \begin{scope}[shift = {(6,-6)}]
    \PythagorasTree{9}{60}
    \node at (.55,-.7) {$\theta=30$};
    \end{scope}
    \end{scope}
    \end{tikzpicture}
\end{center}
Indeed, this puts us precisely in the setup of Definition \ref{maindef}.  For angles $\theta$ in the interval $[30,45]$, we can let $P_\theta\subseteq\bR^2$ be the Pythagorean tree fractal constructed using the right triangle of minimal angle $\theta$ in the recursive step, considered as a filled subset of the plane.  This gives us an entire deforming family of Pythagorean tree fractals, which one can consider as an animation on a screen \cite{Conroy}.\footnote{See Matthew Conroy's animation as $\theta$ varies from 0 to 90: \url{https://www.madandmoonly.com/doctormatt/webimages/animations/pythagoreanFractal2.gif}.}\\

3D printing gives us a new way to observe the transformation of this fractal as $\theta$ varies.  Indeed, our approach is to follow the framework of Section \ref{Theory}, and to consider the deformation space as a solid static animation in $\bR^3$:
\[\cP = \{(x,\theta):\theta\in[30,45],x\in P_\theta\}\subseteq\bR^3.\]
We can use a computer to model $\cP$ using CAD software (we did this using the OpenSCAD \cite{OpenSCAD}), and create a 3D object which can be 3D printed, held in your hand, and studied from all points of view.
\begin{center}
    \includegraphics[width=.85\textwidth]{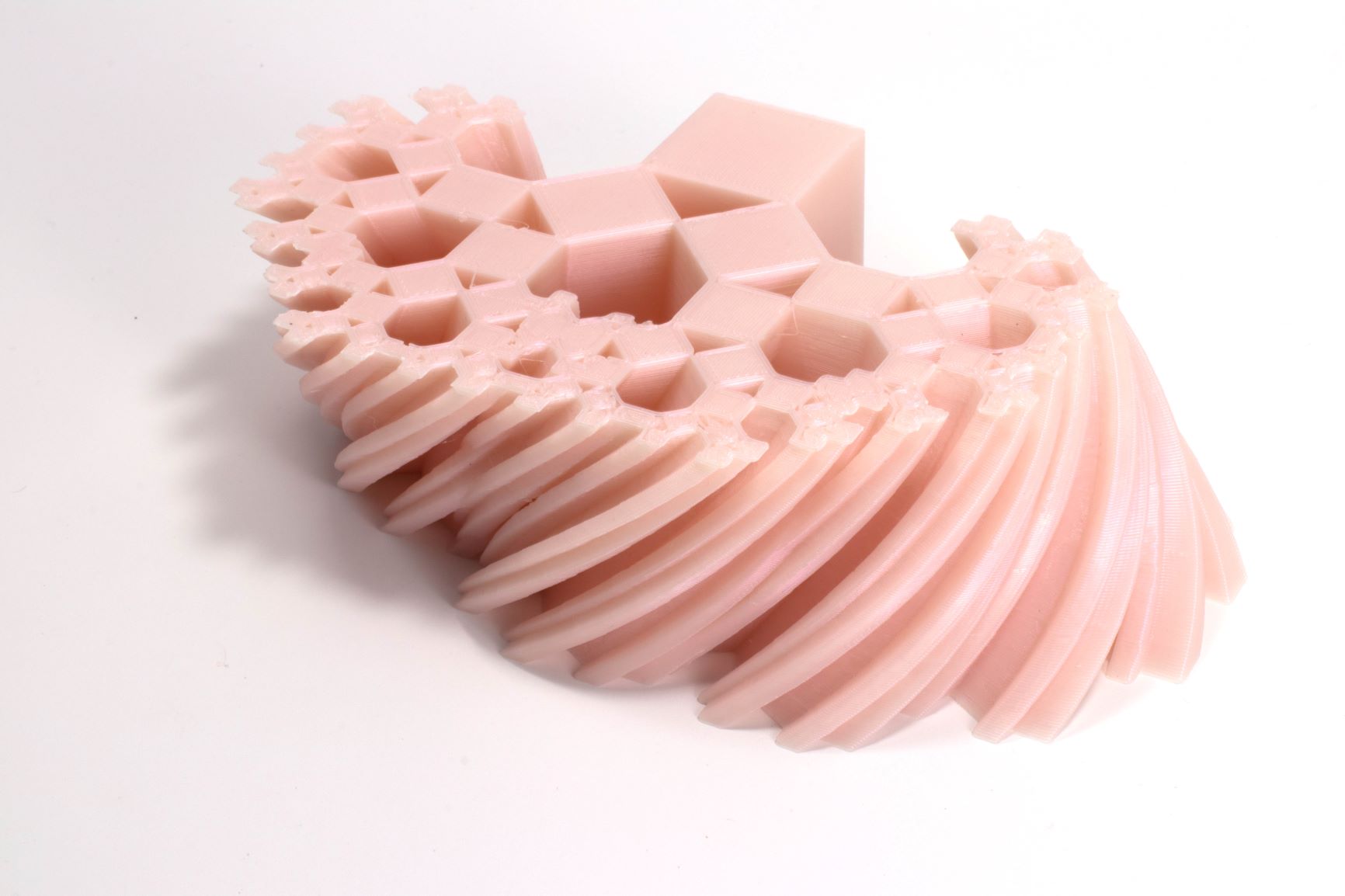}\\
    \vspace{10pt}
    \includegraphics[width=.4125\textwidth]{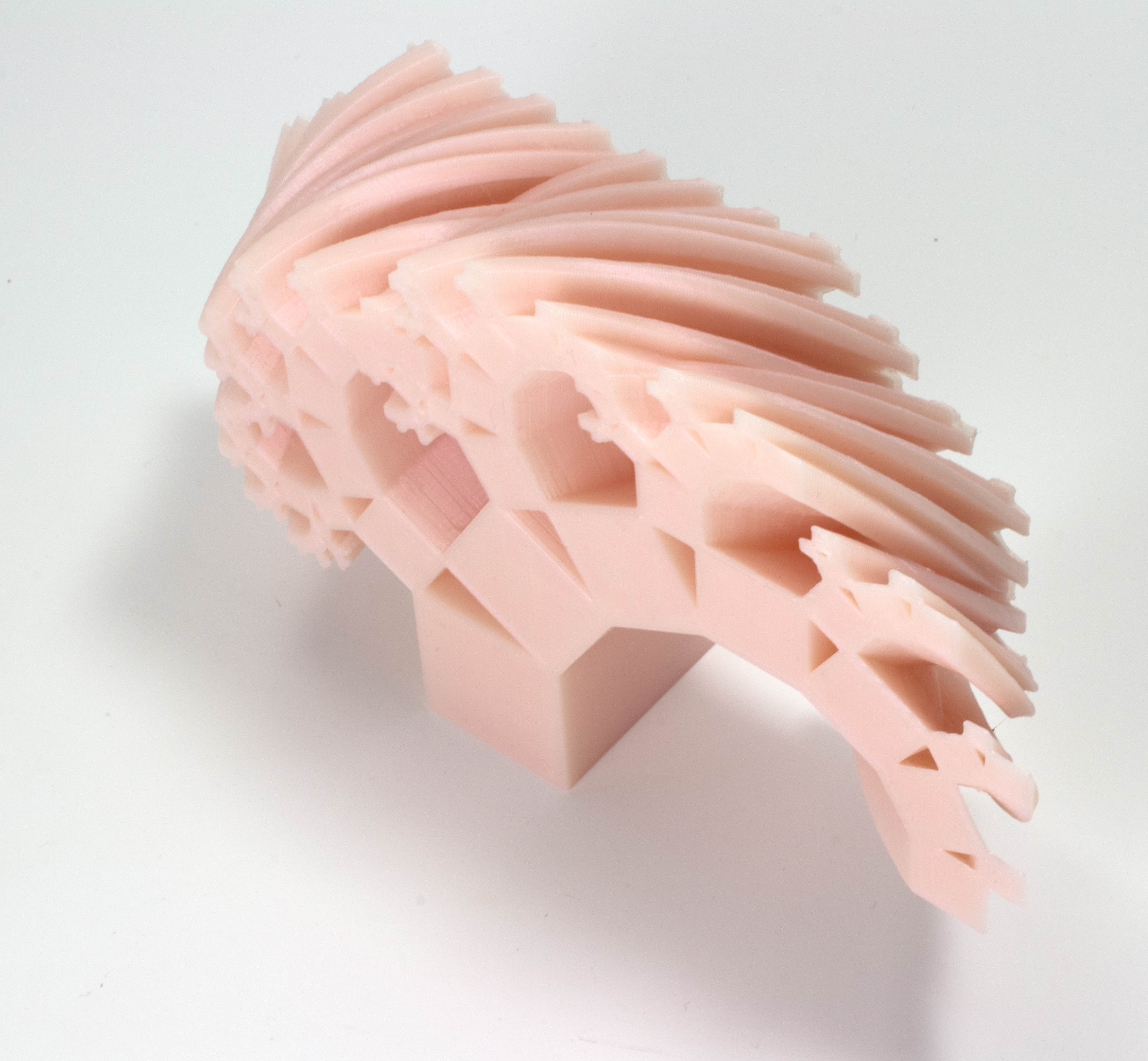}\hspace{10pt}
    \includegraphics[width=.505\textwidth]{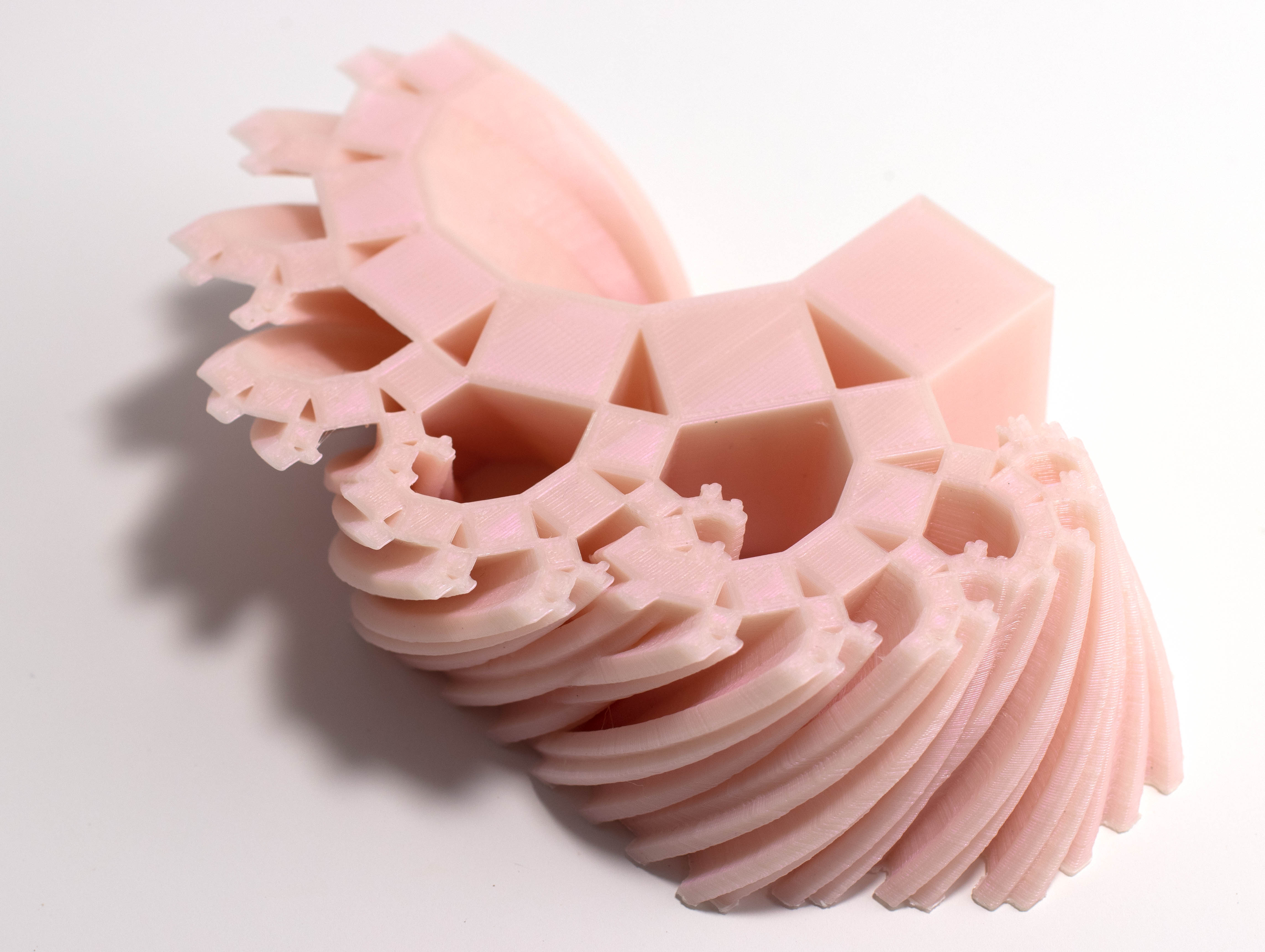}
\end{center}
In this way, 3D printing gives us a new avenue to illustrate perturbations and deformations of this family of fractals.  For one, we get a new dimension where we can observe the movement$-$for example, by considering where the slope becomes more exaggerated, it becomes immediately visible that the change in angle has more of an effect toward the outer edges of the fractal, and more of an effect on the \textit{30 degree} side than the \textit{60 degree side}.  Even more, one obtains a physical and tactile way to experience a mathematical concept, opening new doors for broad mathematical interaction.

\section{Deformations of Julia Sets}\label{juliaSection}
\subsubsection*{Joint with Bernat Espigul\'e}
Let $\bC$ denote the field of complex numbers, fix a complex number $c\in\bC$, and consider the function $f_c:\bC\to\bC$ given by the rule $f_c(z) = z^2+c$.  What happens when we iterate this function?  More explicitly, if we fix some $z_0\in\bC$, what can we say about the following sequence of complex numbers?
\begin{eqnarray*}
z_0&&\\
z_1 &=& f_c(z_0)\\
z_2 &=& f_c(z_1) = f_c(f_c(z_0))\\
z_3 &=& f_c(z_2) = f_c(f_c(f_c(z_0)))\\
\vdots&&
\end{eqnarray*}
A fundamental question asks how the behavior of this sequence changes as we vary $z_0$, or as we vary $c$.  Let's get a basic idea of what can happen by letting $c=0$.
\begin{example}\label{circleExample}
Consider $f_{0}(z) = z^2$.  Then
\[f_0^n(z_0) = \underbrace{f_0\circ f_0\circ\cdots\circ f_0}_{n-\text{times}}(z_0) =
z_0^{2n}.\]
If $|z_0|<1$ then this sequence converges to 0, if $|z_0|>1$ then this sequence escapes to infinity, and if $|z_0|=1$ then this sequence stays on unit circle $|z|=1$.  In particular, the set
\[\{z,f_0(z),f_0^2(z),f_0^3(z),\cdots\}\]
is bounded if and only if $|z|\le1$.
\end{example}
It appears that when iterating $f_c$ from a fixed starting point we can expect the outputs to either escape to infinity, or else to all be contained in a bounded region of $\bC$.  Filled Julia sets consist of the starting points which produce bounded sequences.
\begin{definition}\label{JuliaSetDefinition}
Let $c\in\bC$ be a complex number and consider the function $f_c(z) = z^2+c$.  The \emph{filled Julia set} associated to $c$ is the set:
\[J_c:=\{z\in\bC : \text{ the set }\{z,f_c(z),f_c^2(z),...\}\text{ is bounded}\}\subseteq\bC.\]
\end{definition}
In Example \ref{circleExample} we determined that the filled Julia set $J_0$ is precisely the complex unit disk $\{z:|z|\le1\}\subseteq\bC$.  For other complex numbers, the situation becomes much richer.  For example, below are plots of $J_c$ for various values of $c$\footnote{The 2D images of filled Julia sets in this section were rendered in Mathematica \cite{Mathematica}}.
\begin{center}
    \includegraphics[width=\textwidth]{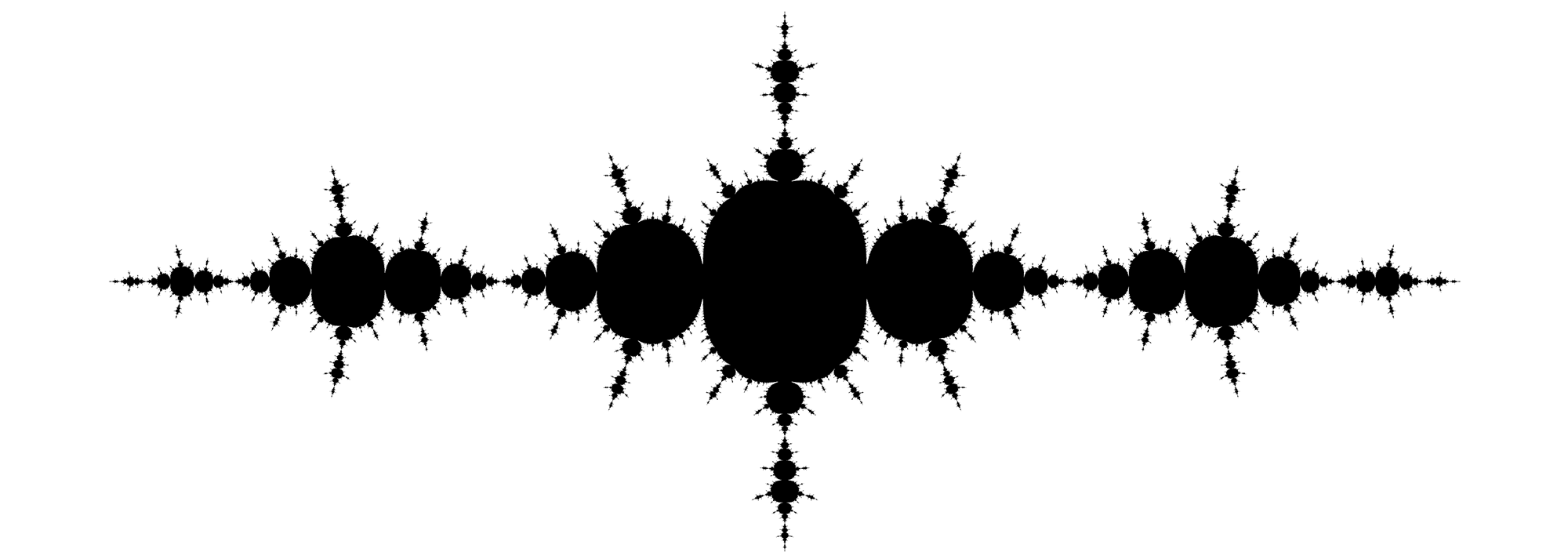}
    $c = -1.25$\\
    \includegraphics[width=.3\textwidth]{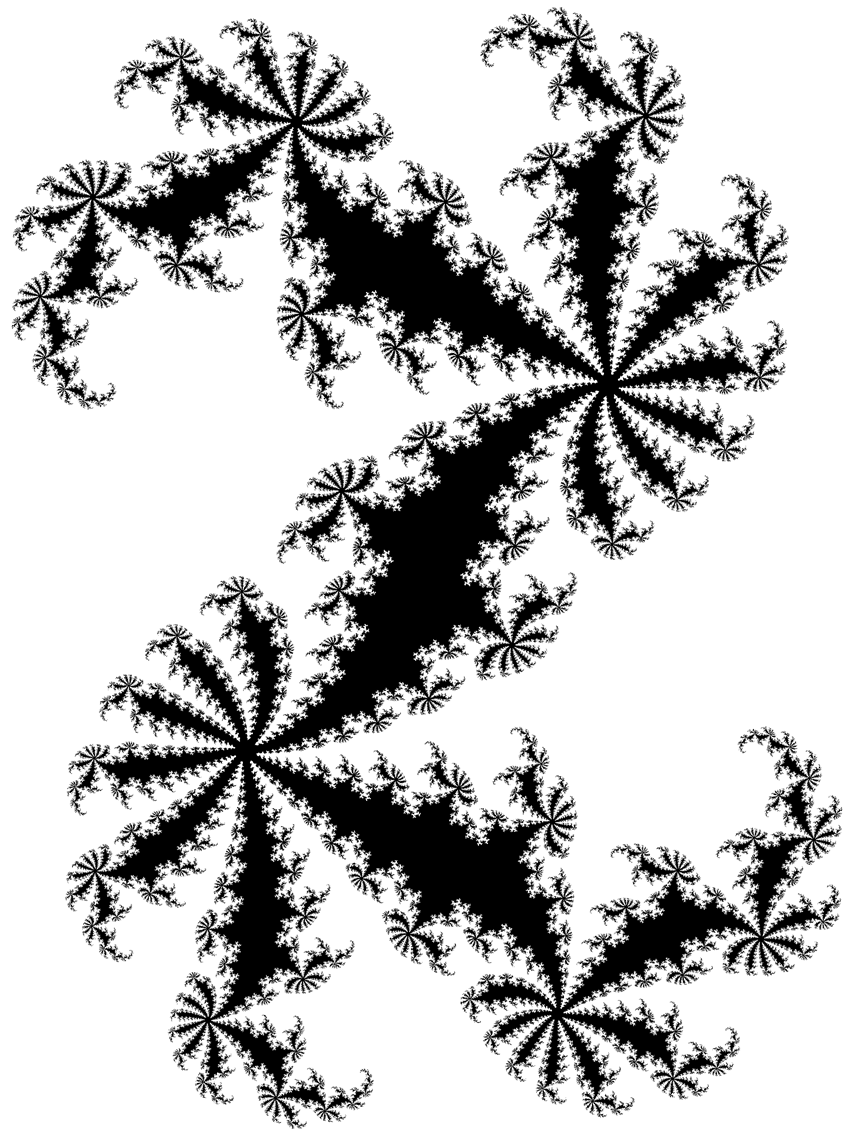}
    \includegraphics[width=.45\textwidth]{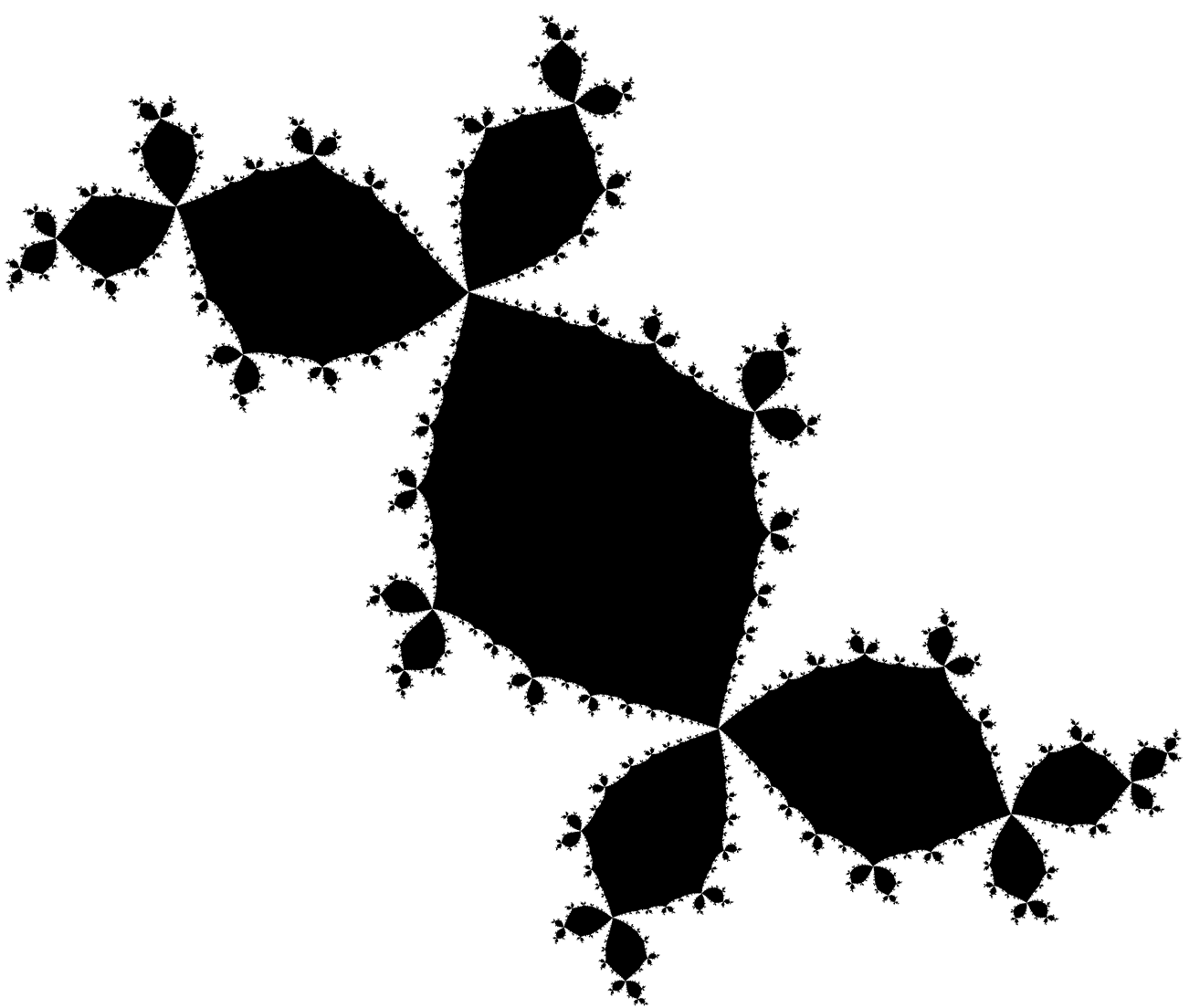}\\
        $c = .365 + .1*i$\hspace{50pt}$c = -0.12+0.74i$
\end{center}
A central question asks how $J_c$ changes as $c$ does.  For example, observe the variation of the filled in Julia sets for polynomials $f_c$ with values $c=-.75$, $c=-.74+.05i$, and $c=-.73+.1i$ respectively.
\begin{center}
    \includegraphics[width=.95\textwidth]{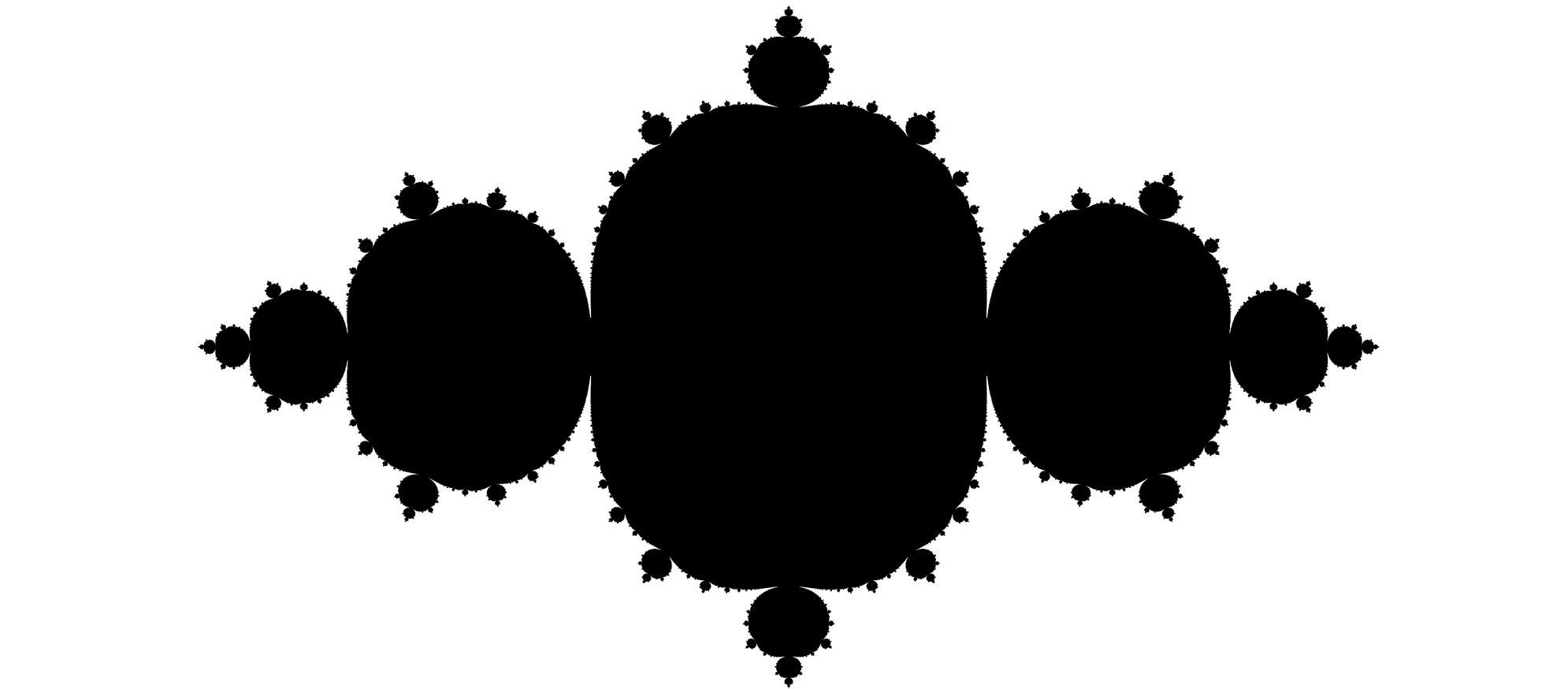}
    \includegraphics[width=.95\textwidth]{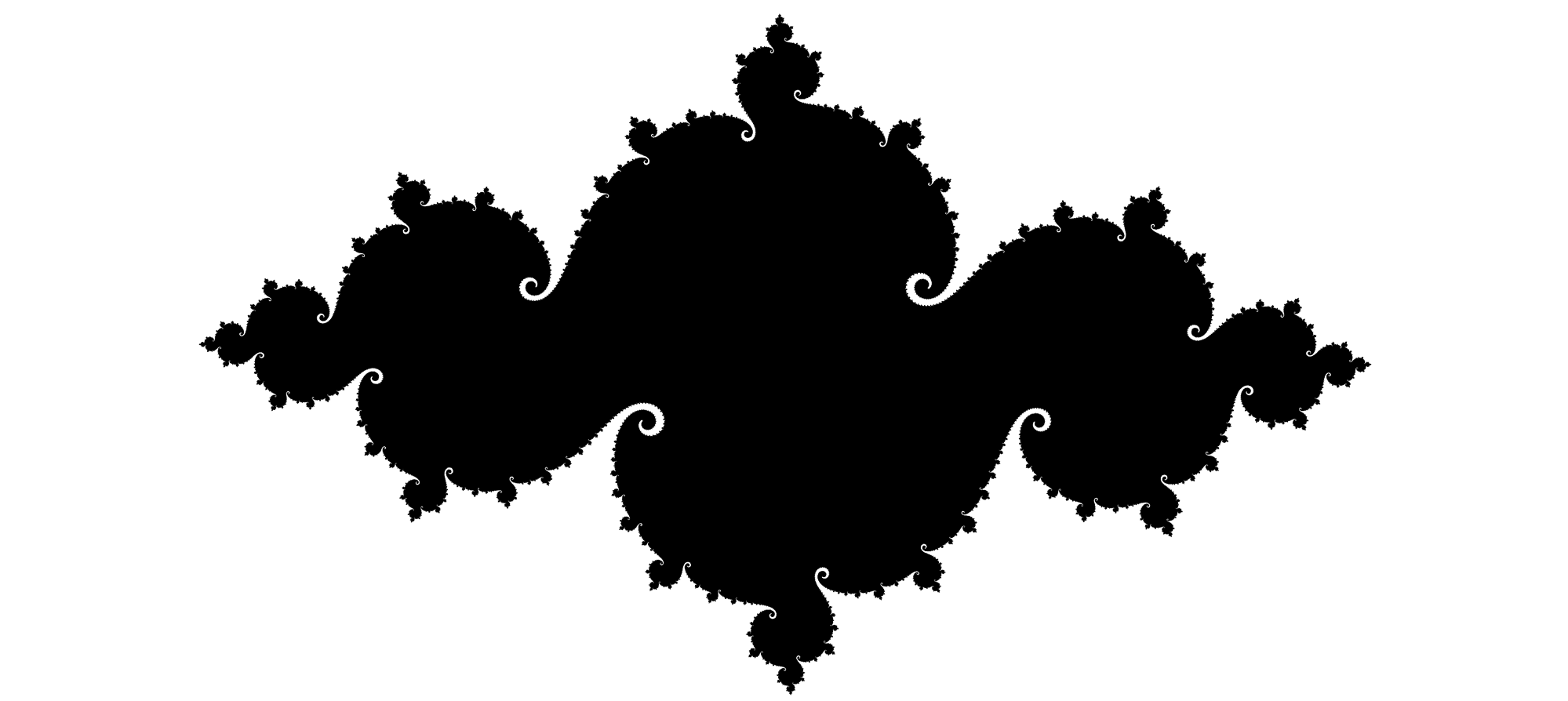}
    \includegraphics[width=.95\textwidth]{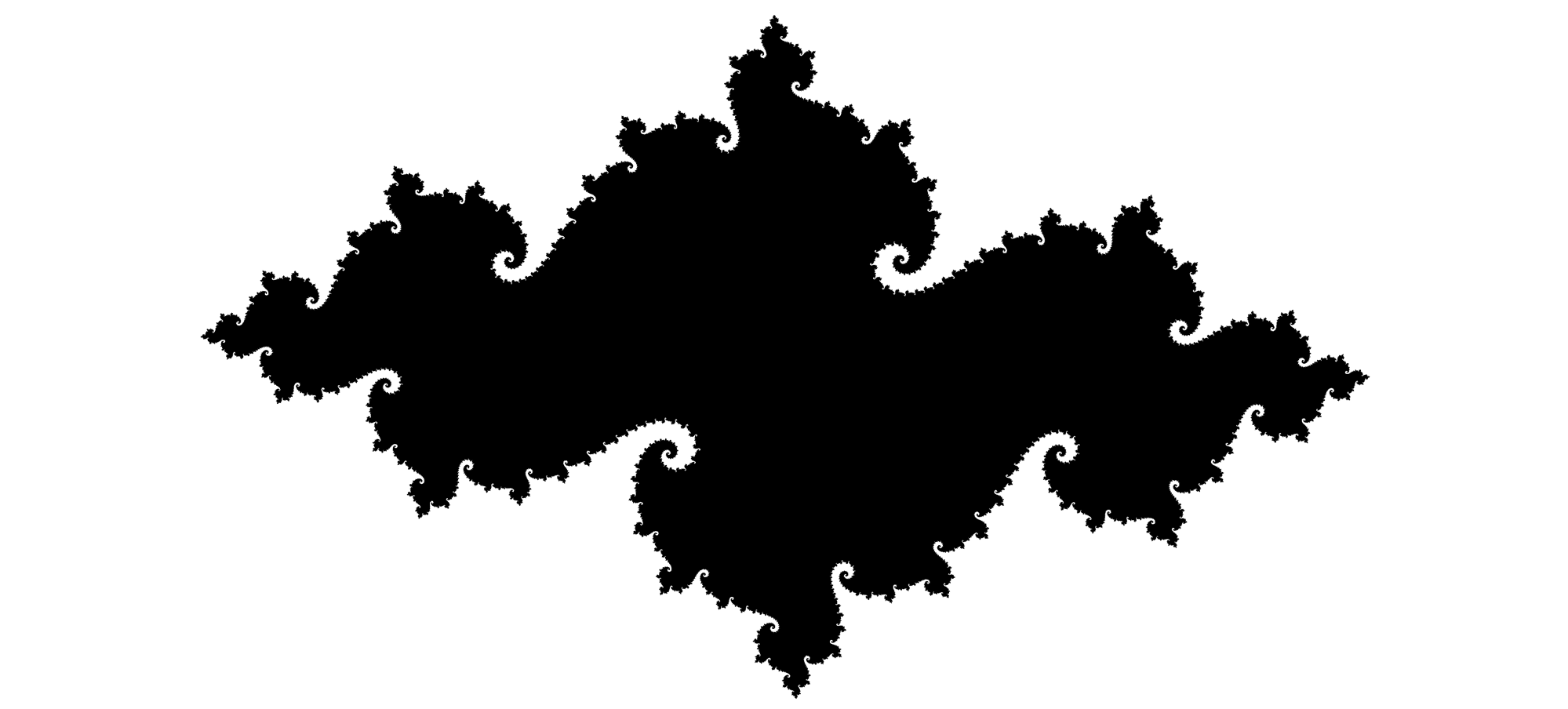}
\end{center}
It appears that if $c$ changes a little bit, then so does $J_c$,\footnote{In fact, the situation is quite a bit more subtle than this.  We discuss this a little further in Section \ref{parabolicImplosions}. The question is studied systematically by Douady \cite{Douady}.} suggesting that we can use these filled Julia sets as the frames of a static animation.  But, as discussed in Section \ref{HDF}, there is an obstacle.   This is because the varying parameter $c$ is a complex number, and so we obtain a 2-dimensional stack of frames.  When we form a deformation space exactly as in Definition \ref{maindef}, we obtain a subset of $\bC\times\bC$:
\[\cJ = \{(c,z):c\in\bC,z\in J_c\}\subseteq\bC^2.\]
As a hypersurface in 4D space, this object is impossible to 3D print.  In order to get something we \textit{can} print, we will choose a 1-dimensional path through the parameter space (c.f. Section \ref{HDF}).  But which paths should we choose?  To answer this it is important to understand a little more about the geometry of the parameter space of filled Julia sets.\\

Each filled Julia set we have looked at so far has been connected and has a clearly defined interior, but for many values of $c$ the set $J_c$ is a totally disconnected Cantor set \cite[Chapter 9]{Milnor}.  For the purposes of 3D printing it makes sense to focus on values $c$ where $J_c$ is connected, which is a fascinating and celebrated collection of complex numbers in its own right.
\begin{definition}\label{Mandlebrot}
The \textit{Mandelbrot set} is the set:
\[\cM:=\{c\in\bC: J_c\text{ is connected}\}\subseteq\bC\]
\end{definition}
A more familiar definition of the Mandlebrot set is perhaps
\[\{c\in\bC: \{0,f_c(0),f_c^2(0),...\}\text{ is bounded}\} = \{c\in\bC: 0\in J_c\}.\]
The fact that these are equivalent is \cite[Theorem 9.5, Appendix G]{Milnor}.  The Mandlebrot set famously has the following shape:
\begin{center}
    \includegraphics[width=.75\textwidth]{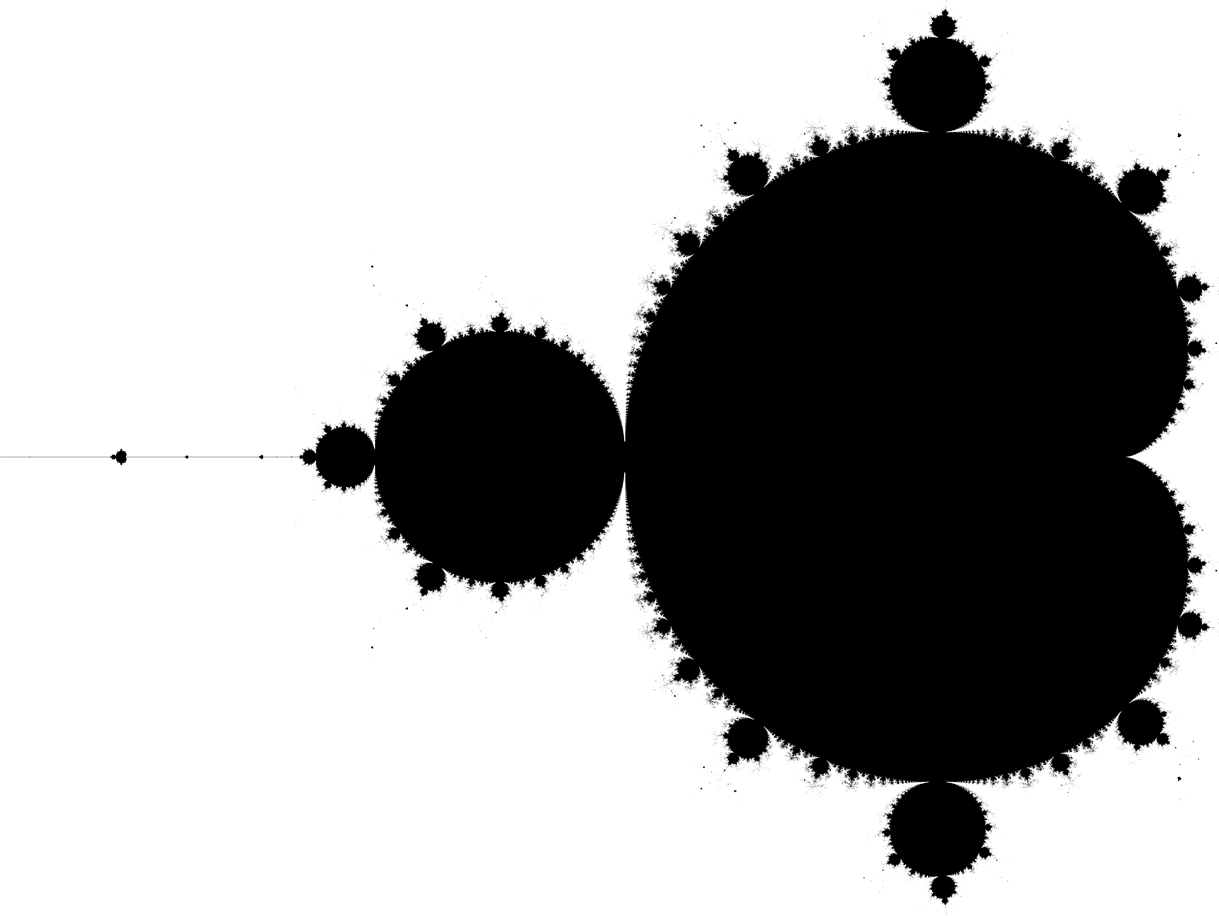}
\end{center}
To cut down our deformation space $\cJ$ by one dimension and obtain something 3D printable, we mirror Section \ref{HDF} and pick a 1-dimensional subset of the Mandlebrot set defined by a path.  As the behavior of the Julia set $J_c$ is closely related to the location of $c$ in the Mandlebrot set, we will use the geometry of the Mandlebrot set to inform our choice of path.
\subsection{Turning Paths into Prints}\label{P2P}
Let $I\subseteq\bR$ be an interval, and consider a continuous map
\[\gamma:I\to\bC,\]
whose image is contained in $\cM$.  For each $t\in I$, we can consider the filled Julia set $J_{\gamma(t)}$.  This gives us a family of filled Julia sets parametrized by $t$, putting us squarely in the situation of Definition \ref{maindef}.  Therefore the deformation space associated to this data is a static animation:
\[\cJ_\gamma:=\{(z,t):t\in I, z\in J_{\gamma(t)}\}\subseteq\bC\times\bR\approx\bR^3.\]
This sets up a general framework for creating static animations illustrating the deformation of filled Julia sets associated to a path in the Mandlebrot set.  It only remains to choose a path.  The boundary of the Mandlebrot set partitions the collection of filled Julia sets into two types: the filled Julia sets that are connected, and those that are not.  In particular, it behaves like a critical point, and we begin by studying how filled Julia sets deform near the boundary, \textit{or along it}.  Our initial exploration was to traverse the upper boundary of the main cardioid.
\begin{center}
\includegraphics[width=.6\textwidth]{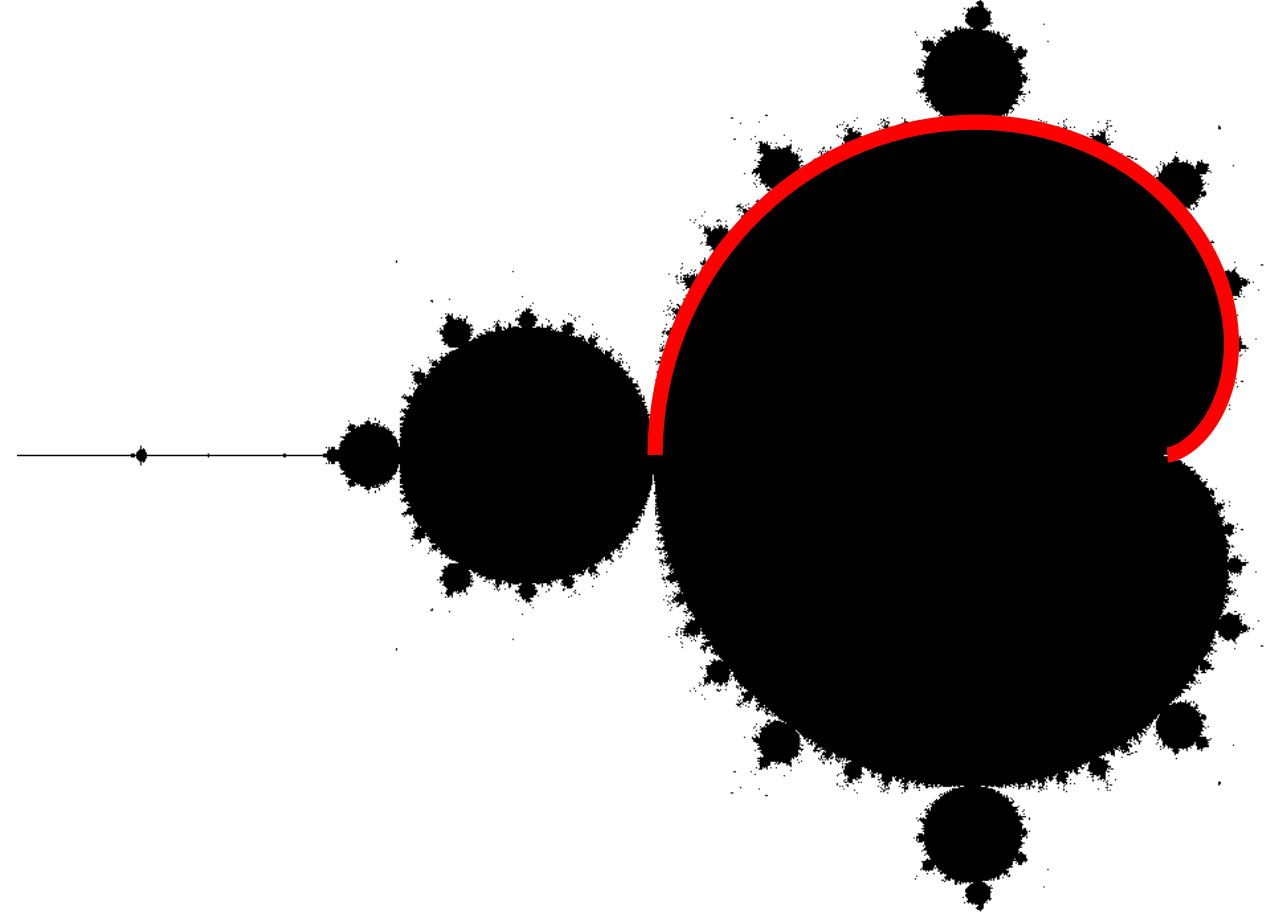}
\end{center}
One can deduce from \cite[Chapter 10]{Milnor} that the boundary of the main cardioid can be parametrized as follows.\begin{eqnarray*}
    x(t) &=& \frac{1}{2}\cos(t)(1-\cos(t)) + \frac{1}{4},\\
    y(t) &=& \frac{1}{2}\sin(t)(1-\cos(t)).
\end{eqnarray*}
Then the path $\gamma:[0,\pi]\to\bC$ given by $\gamma(t) = x(t) + iy(t)$ parametrizes the upper boundary of the main cardioid.  Associated to this path there is a static aimation $J_\gamma\subseteq\bC\times\bR$, a 3D print of which is pictured below.
\begin{center}
    \includegraphics[width=.4\textwidth]{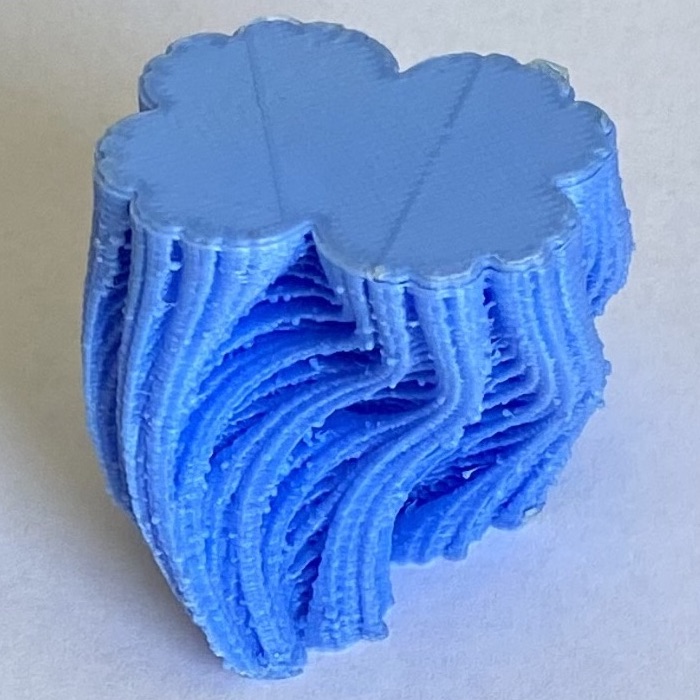}\hspace{10pt}
    \includegraphics[width=.4\textwidth]{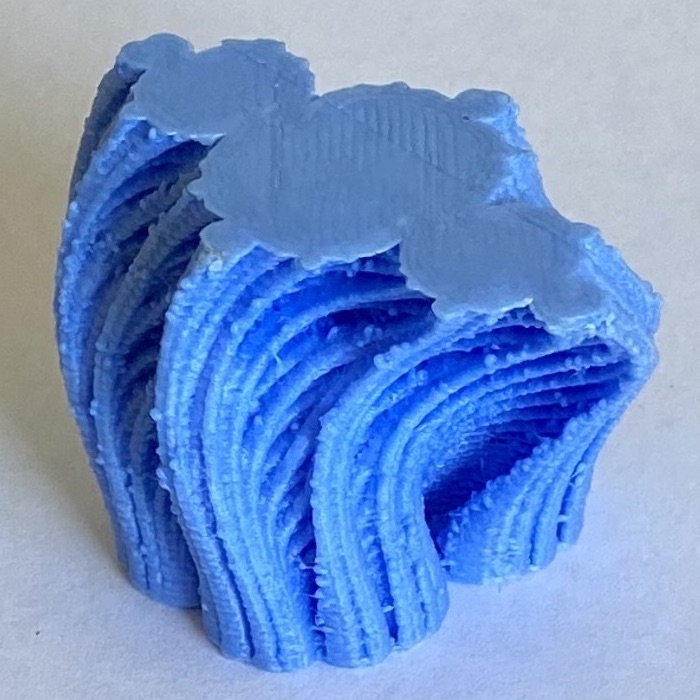}
\end{center}
\subsection{Parabolic Implosions}\label{parabolicImplosions}
We suggested in Section \ref{Theory} that we were mostly interested in static animations with continuously deforming frames, but is this true for Julia sets?  This question was studied by Douady \cite{Douady} who showed that in general, it is not.  The fact is, there are some special points at which $c\mapsto J_c$ does not vary continuously in $c$, even along our path $\gamma$ traversing the main cardioid.  To illustrate an example of this we look at the filled Julia sets associated to $\gamma(0)$ and $\gamma(\varepsilon)$ for some very small $\varepsilon$ (in this case $\varepsilon = .075$).  These are the frames $J_c$ associated with the complex numbers $c_0 = .25$ and $c_{\varepsilon} =0.251402 + 0.000105321 i$ respectively.
\begin{center}
    \includegraphics[width=.25\textwidth]{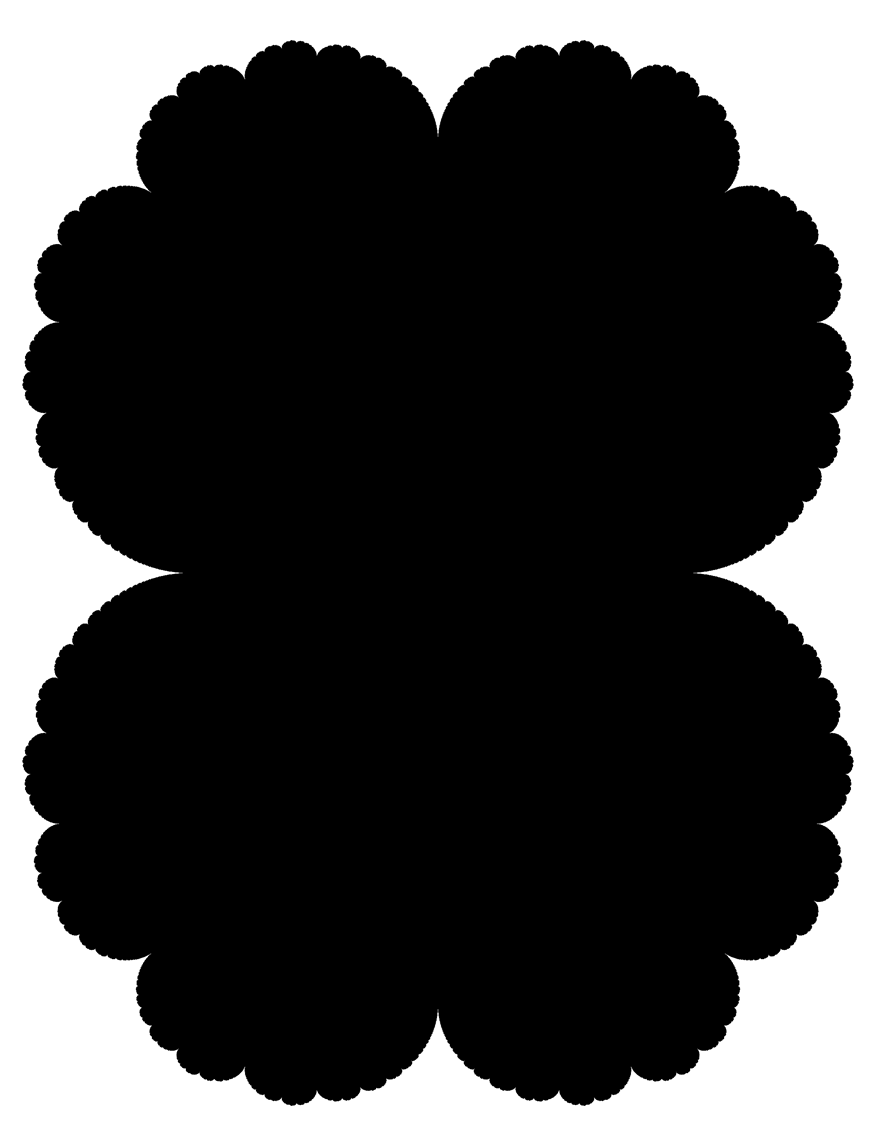}\hspace{100pt}\includegraphics[width=.25\textwidth]{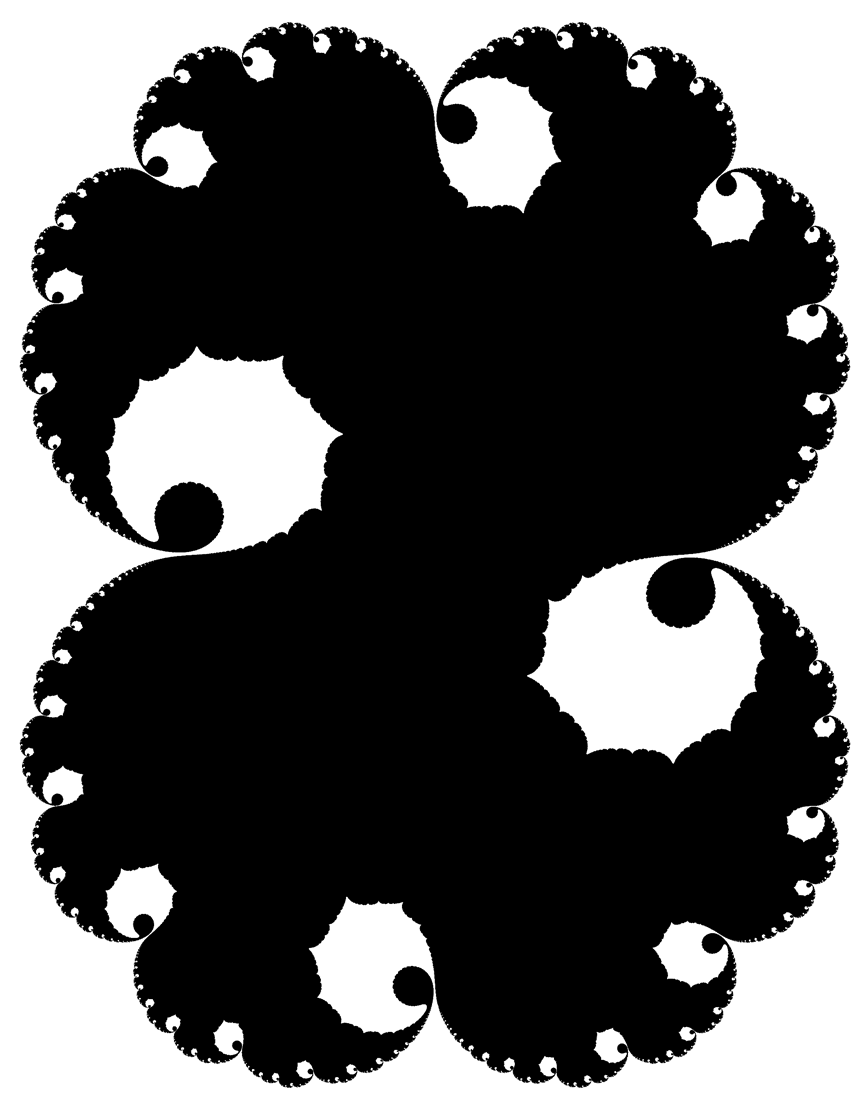}
\end{center}
Although the overall shape is similar, there is a sudden appearance of a number of gaps.  This is an example of a phenomenon called \textit{parabolic implosion} \cite[Section II]{Douady}, which is of considerable interest \cite{Parabolic}, but this curious phenomenon isn't clearly visible in our model, as the filled top layer of the print obstructs the view of the suddenly appearing gaps.  One way to see the parabolic implosions is to instead consider a restriction of our path to a slightly smaller interval.  There is a parabolic implosion also at $\gamma(\pi)$, so we consider the path $\gamma:[\varepsilon,\pi-\varepsilon]$ for very small $\varepsilon$.   The static animation associated to this path allows a better view of formation and deformation of the parabolic implosions, and produces a striking visual effect.
\begin{center}
    \includegraphics[width=.67\textwidth]{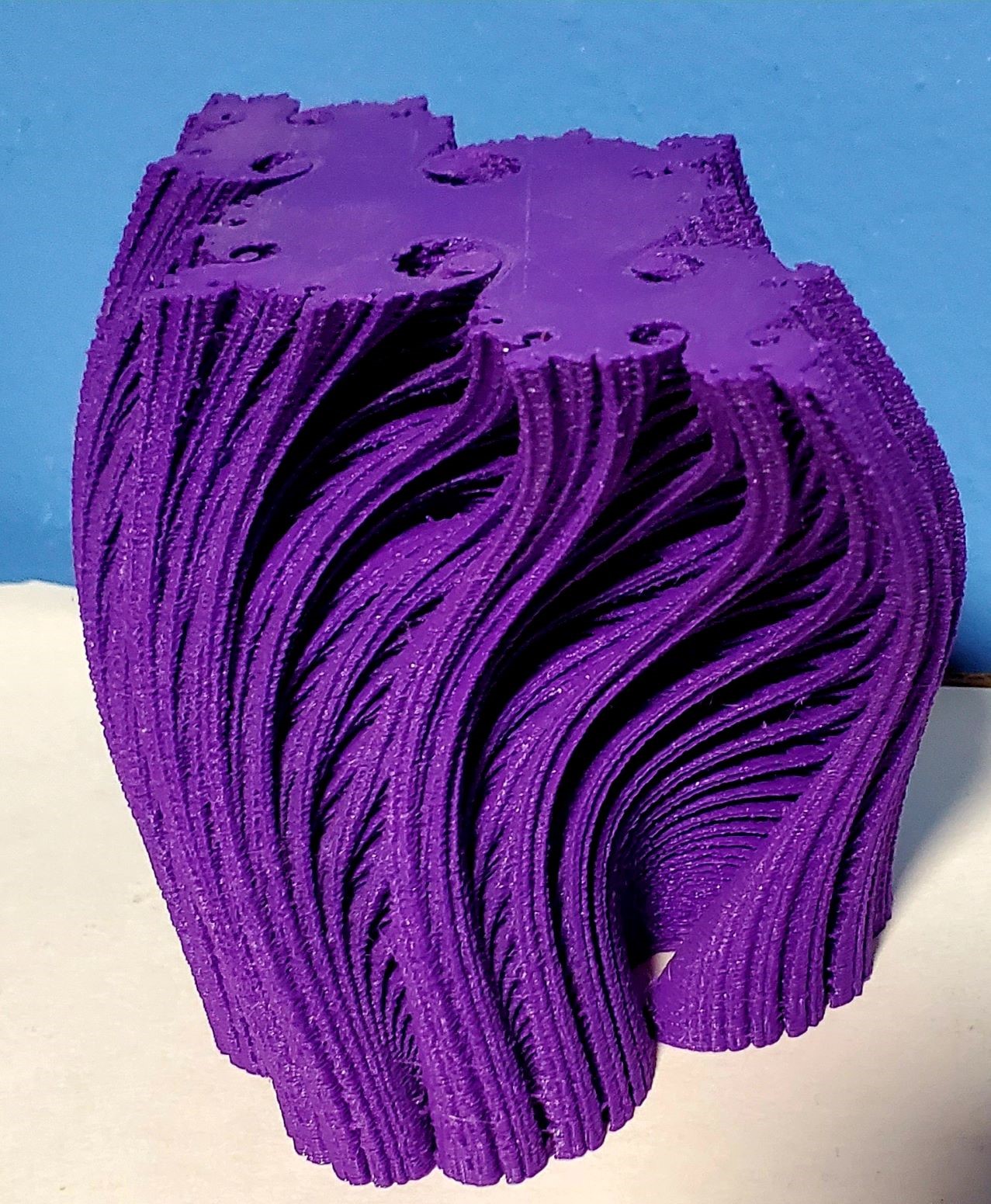}\\
    \vspace{5pt}
    \includegraphics[width=.4\textwidth]{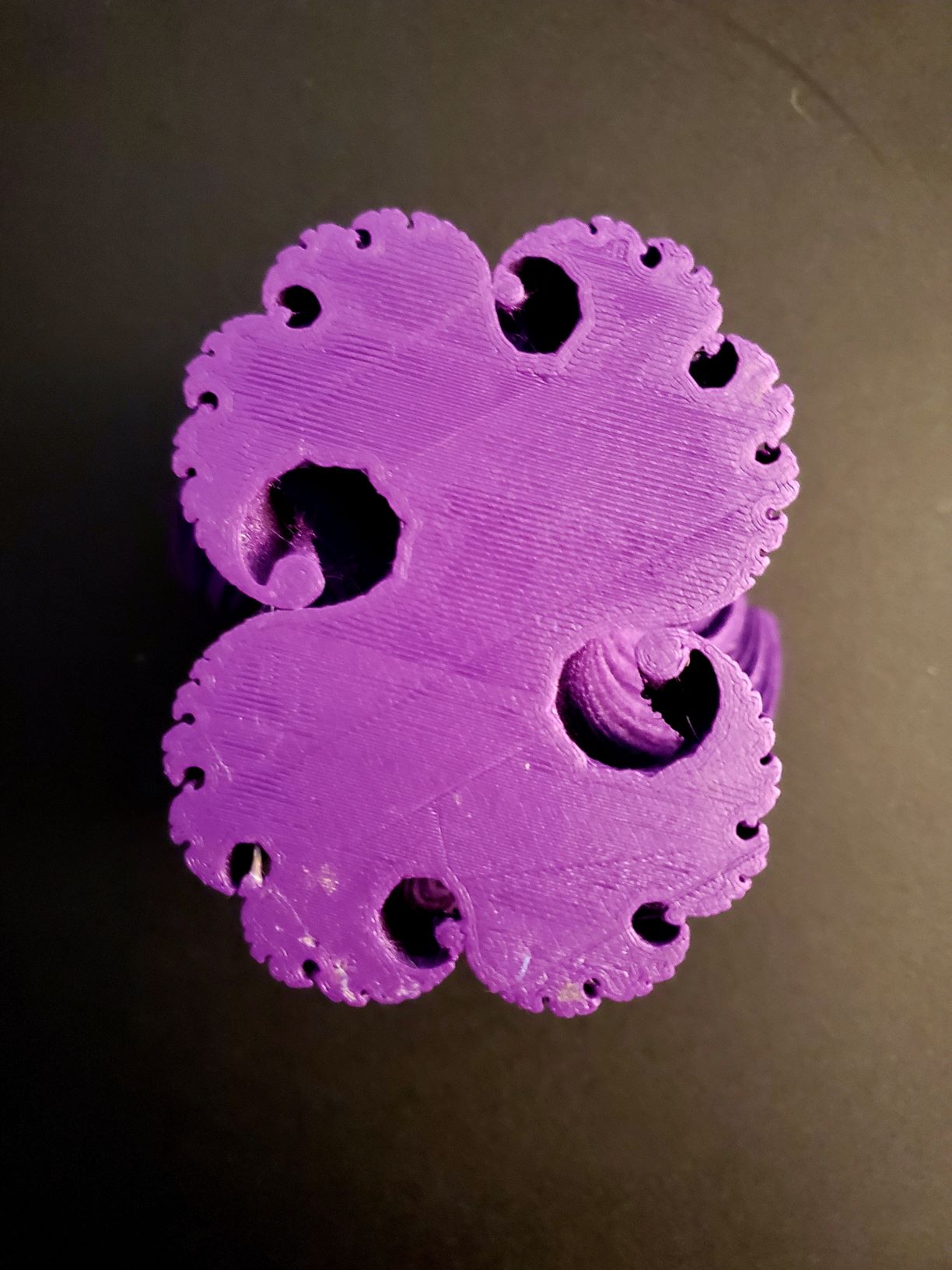}
    \includegraphics[width=.4475\textwidth]{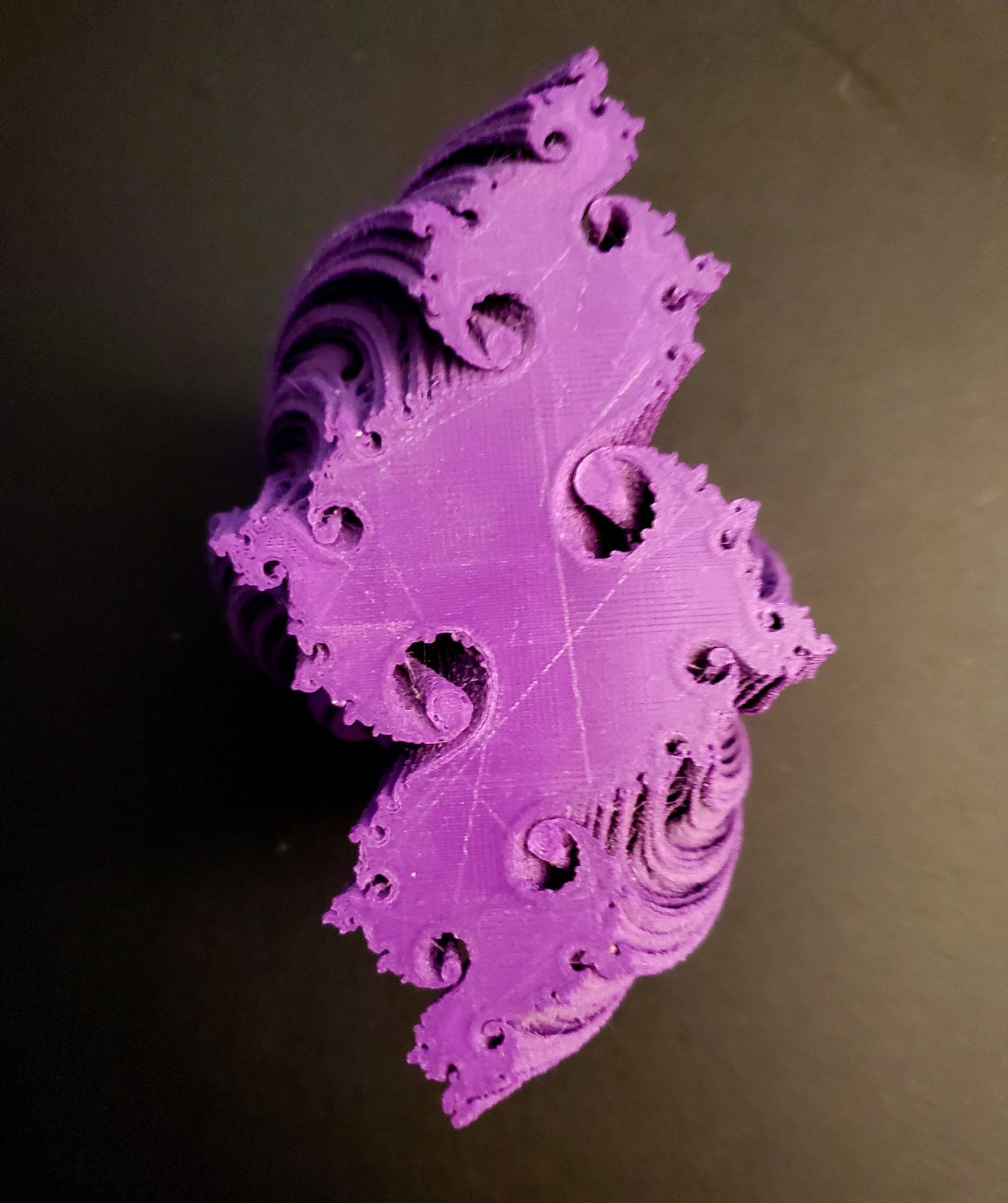}
\end{center}
Having a 2-dimensional parameter space initially seemed like an obstacle in creating static animations, but the philosophy of Section \ref{HDF} turns this into a richness instead.  The freedom of having a static animation associated to each path allows us to focus on interesting areas and study them from multiple points of view.  We exploited this by choosing paths that remain in the vicinity of these parabolic points, allowing us to do local explorations around the parabolic implosions of the cusp ($c\approx\frac{1}{4}$).
\begin{center}
    \includegraphics[width=.45\textwidth]{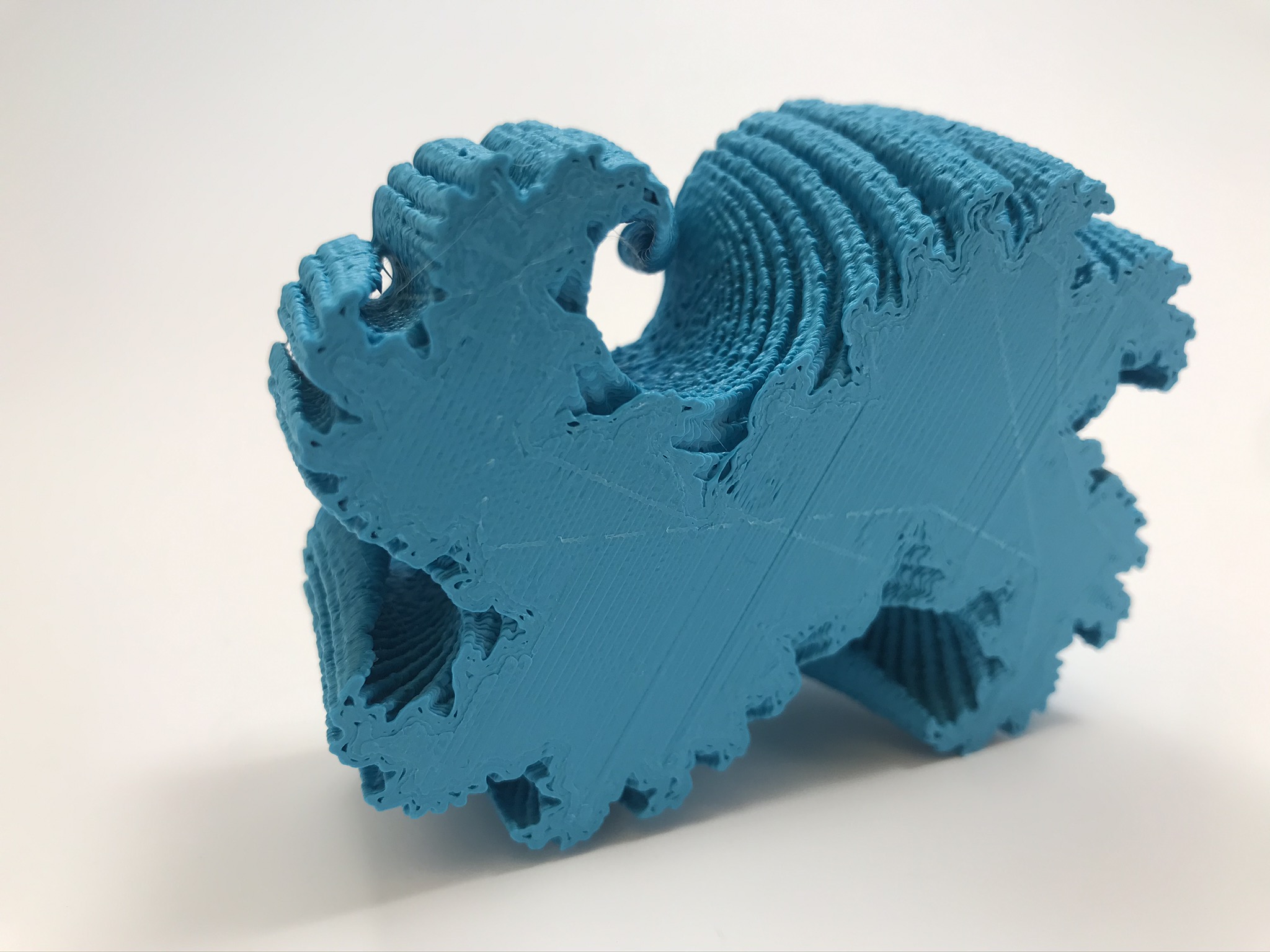}
    \includegraphics[width=.45\textwidth]{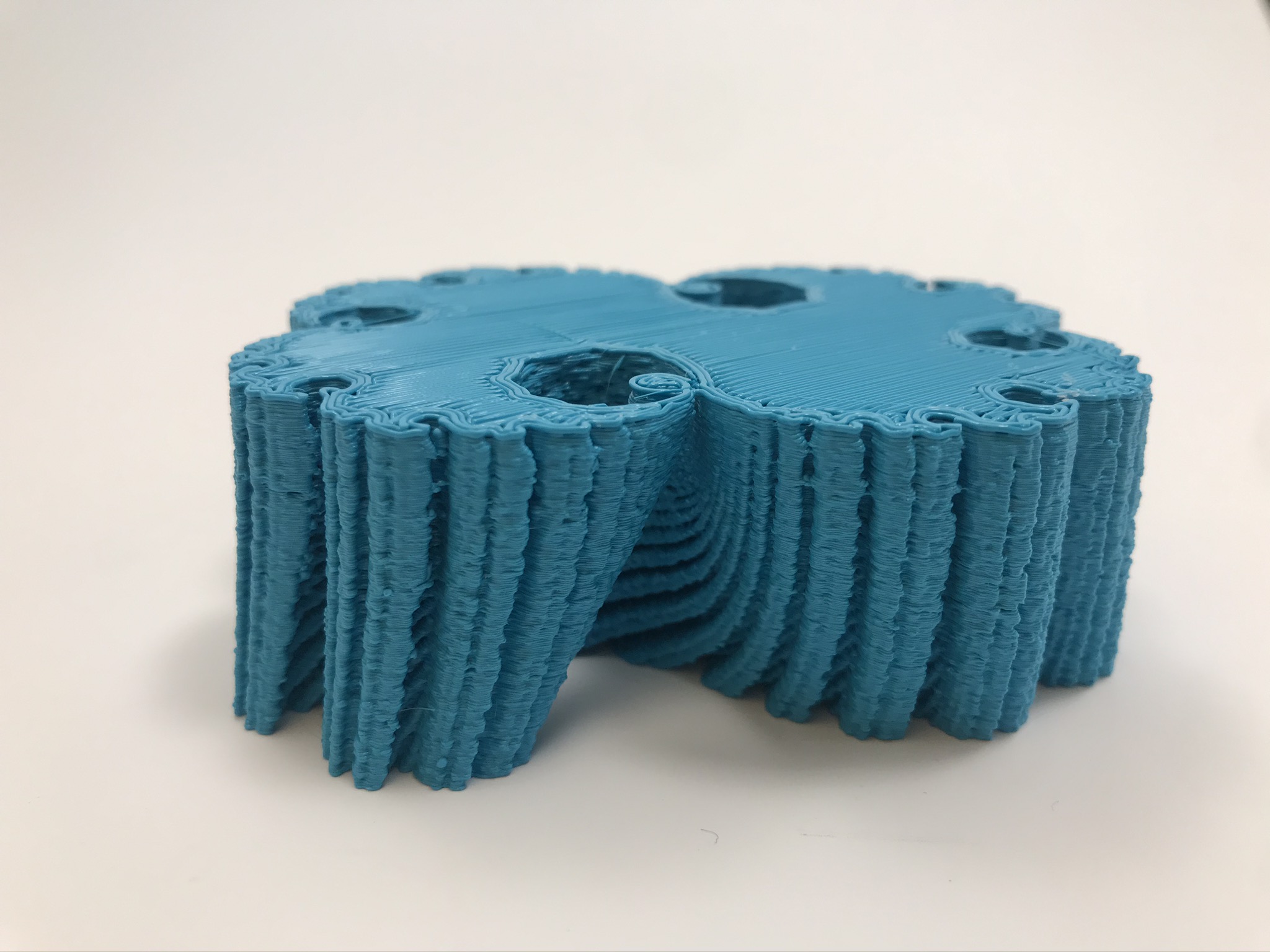}
\end{center}
We do a similar exploration around the parabolic implosion near the so-called Basilica ($c\approx-\frac{3}{4}$).
\begin{center}
    \includegraphics[width=.45\textwidth]{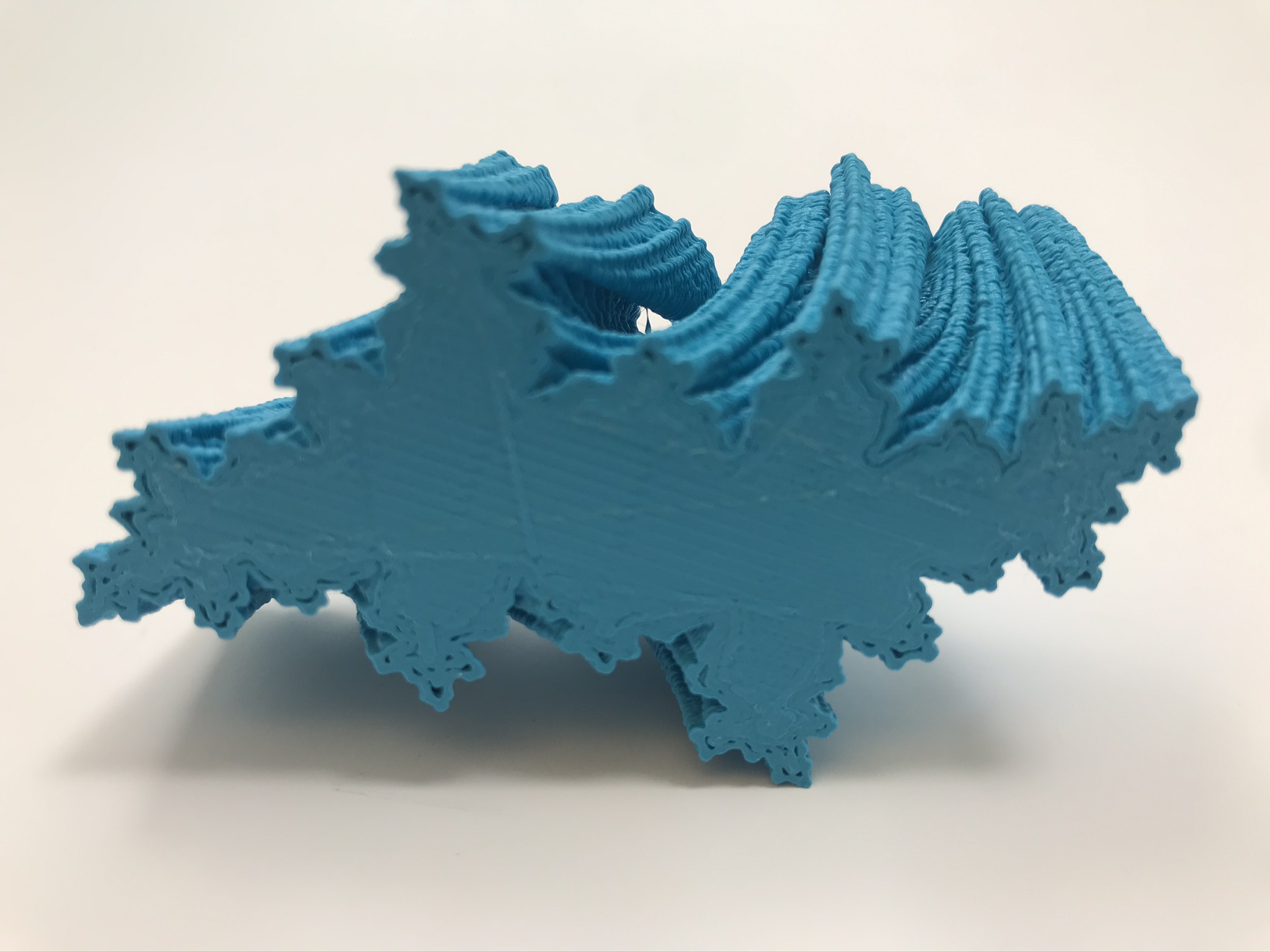}
    \includegraphics[width=.45\textwidth]{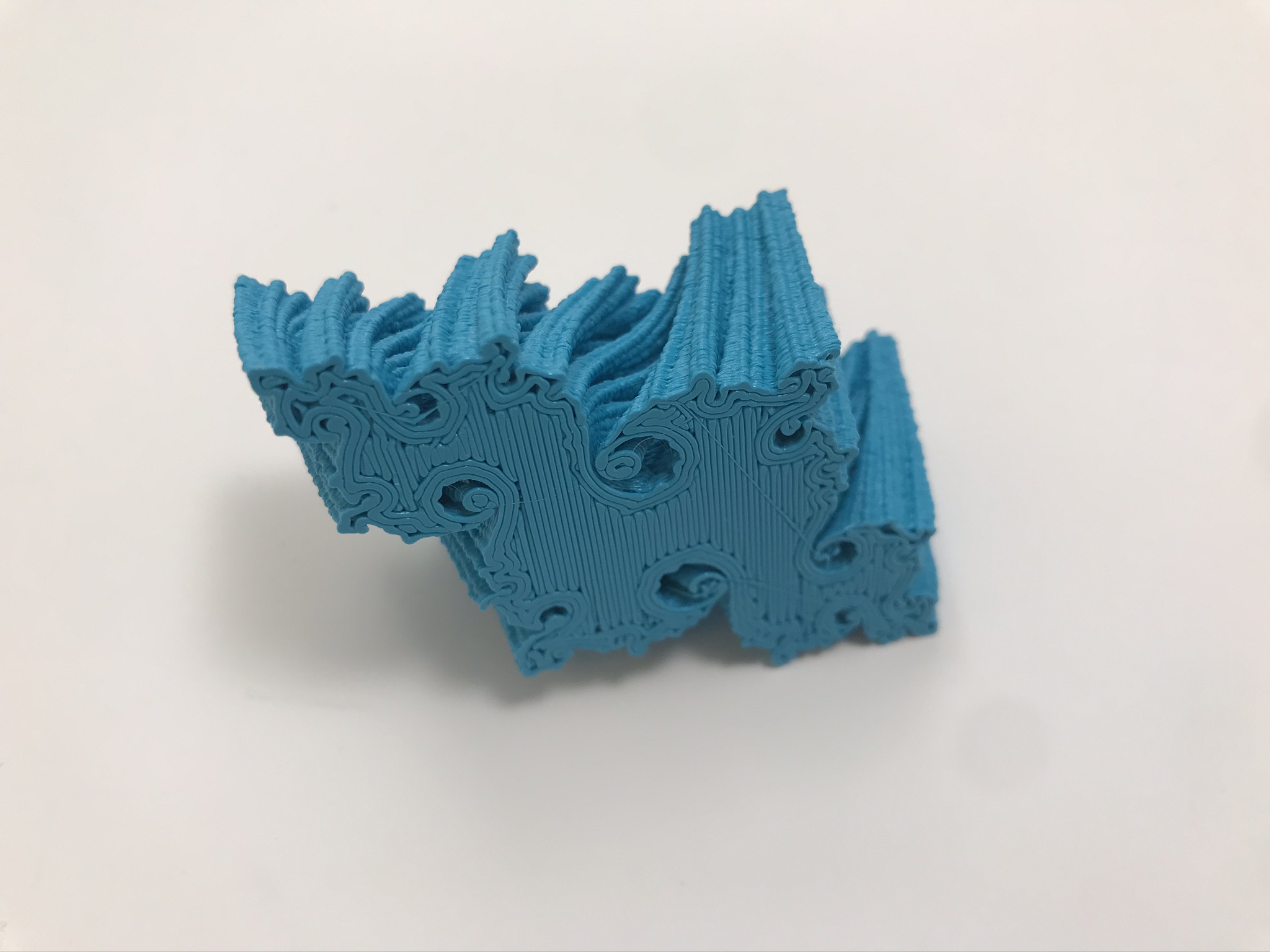}
\end{center}
\subsection{Non-Connected Julia Sets}
So far we have restricted ourselves to considering filled Julia sets $J_c$ for values of $c$ in the Mandlebrot set.  What happens if we try to also consider $J_c$ for $c\in\bC\setminus\cM$?  These Julia sets will not be connected, and in fact, the Julia-Fatou theorem says that $J_c$ is either connected, or else totally disconnected and homeomorphic to a Cantor set \cite[Theorem 9.5]{Milnor}.  We experimented with how our methods (cf. Section \ref{P2P}) adapted to a path which exited the Mandlebrot set, and obtained a curious result.
\begin{center}
    \includegraphics[width=\textwidth]{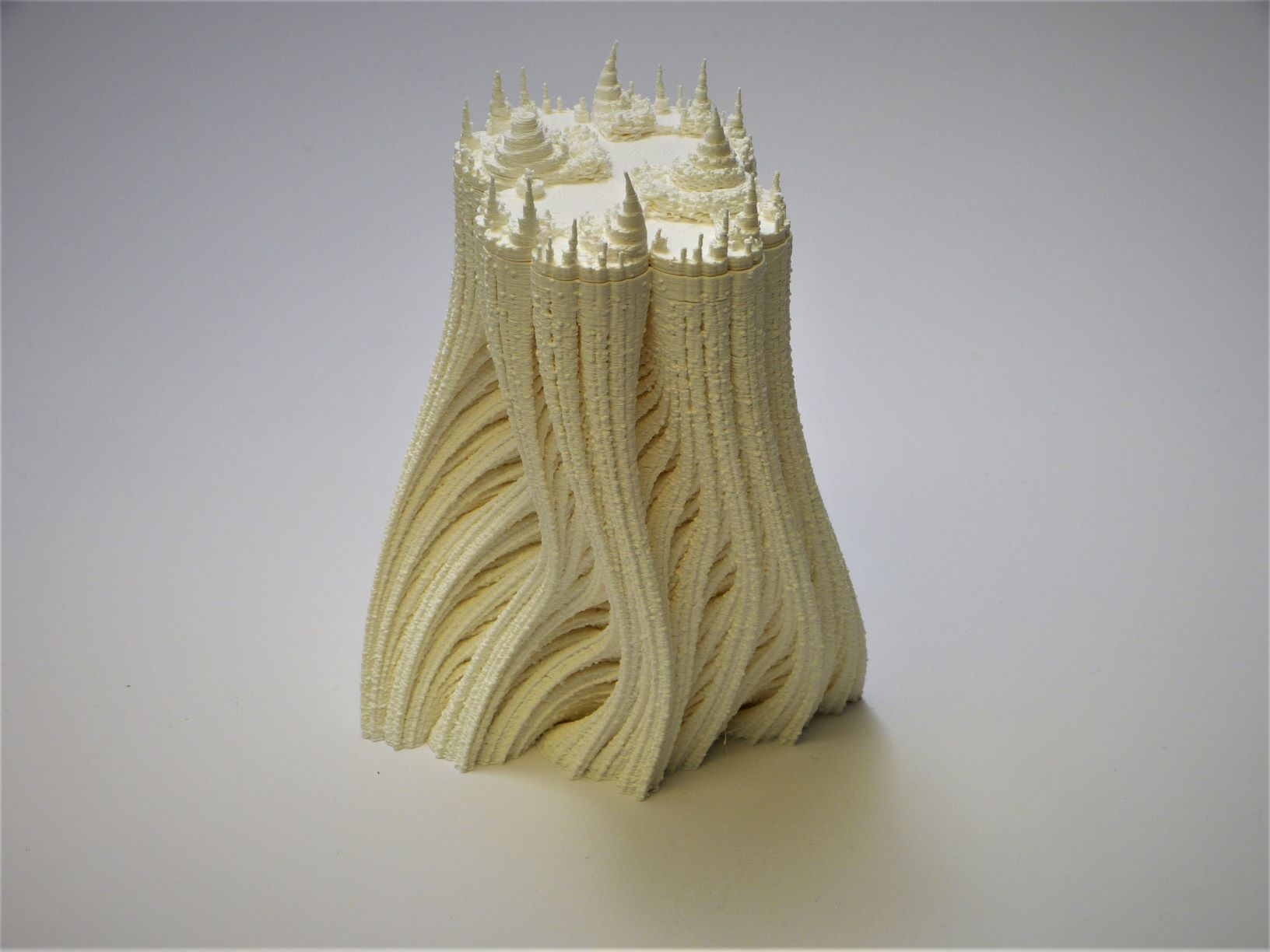}
\end{center}
Where do these towering mountains come from?  In some sense, they are an artifact of computer graphics.  What happens is that our method for generating images of Julia sets uses density to color pixels.  That is, if a pixel contains many points in the Julia set, it will be colored black, and otherwise it will be left white.  As a result, for complex numbers $c$ which are very close to the boundary of $\cM$, the screen displays $J_c$ as a union of several connected components, since these are where there are high densities of points in the Julia set.  For example, here is our image for $c=.28+.01i$, a complex number outside of the Mandlebrot set, but just barely.
\begin{center}
\includegraphics[width=.5\textwidth]{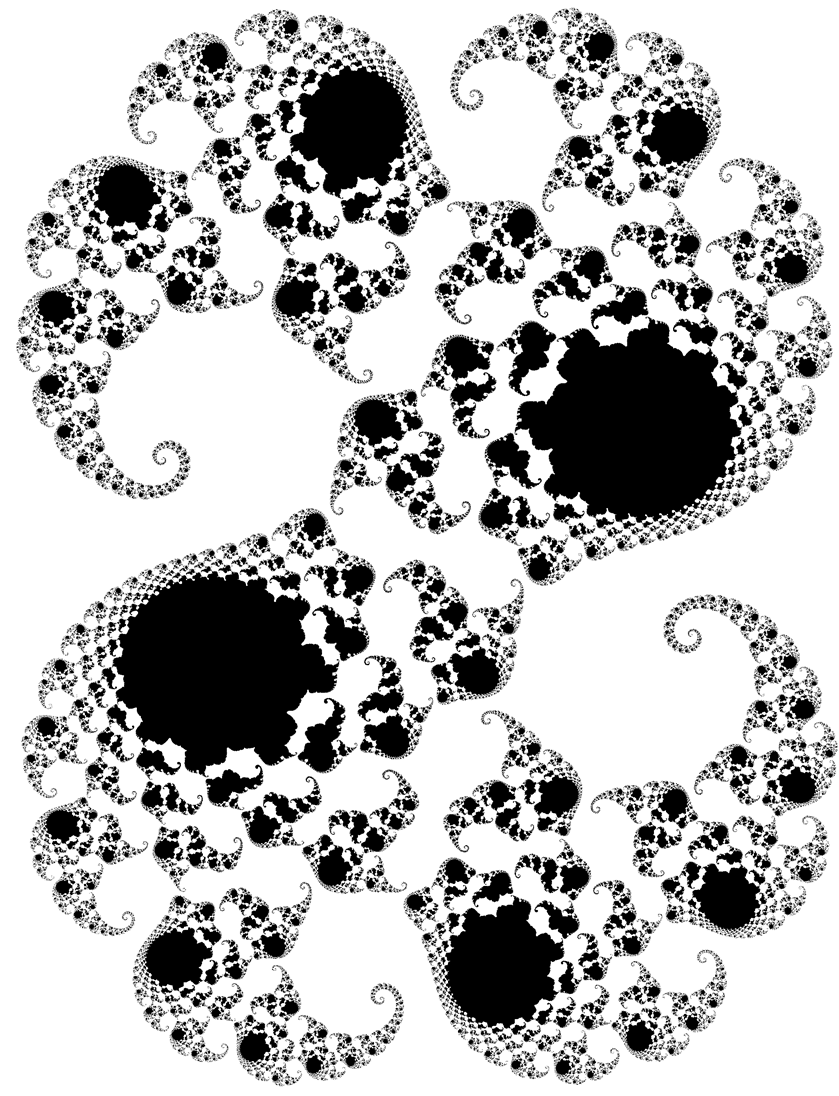}
\end{center}
\textit{We emphasize that this is an artifact of computer graphics}, the set $J_c$ is totally disconnected.  What appears to be a connected region is merely a part of the plane where there is a high concentration of points which are bounded under the dynamical system of iterating $f_c$.  As $c$ progresses further away from $\cM$, these high concentration areas spread out and we can observe a rapid progression toward a more sparse Julia set.  It is worth remarking that this print, although not mathematically rigorous in the purest sense, provides mathematical insights about concentrations of points in Cantor-like Julia sets and how they deform and spread out\\

The fact that an artifact of computer graphics can both create interesting aesthetics and contribute meaningfully to our understanding is interesting in its own right, although beyond the scope of this article.  Discussions about techniques, artifacts, and limitations of generating Julia sets using computer graphics are discussed by Milnor \cite[Appendix H]{Milnor}, and in great depth by Ch\'eritat \cite{Cheritat}.
\subsection{Other Directions}\label{juliaMore}
The static animations of filled Julia sets exhibited so far have been limited to paths on and near the boundary of the main cardioid,  but there is an entire menagerie of static animations that can be constructed simply by choosing paths through the Mandlebrot set.  We invite the reader to explore these possibilities, and to this end we have included a Mathematica notebook called \verb|JuliaPlots.nb| in the public github repository accompanying this paper \cite{github}.   This notebook allows the user to specify a parametrized path through $\bC$, and will generate the frames of the associated static animation.  The reader can then follow the procedure described in Section \ref{How?} to construct the resulting static animation. For now, we suggest a few of the many paths worth exploring, and highlight a related project led by Caroline Davis \cite{log(iu)}.

\subsubsection{Traversing Other Boundaries}
To get a more complete picture of how filled Julia sets deform along the boundary of the Mandlebrot set we should also consider static animations associated to the boundaries of other bulbs.  For example, the bulb directly left of the main cardioid (called the \textit{period 2 hyperbolic component}) is a perfect circle of radius 1/4 centered at $-1$.  Below is the static animation associated to traversing the upper semicircle of this component.
    \begin{center}
        \includegraphics[width=.3\paperwidth]{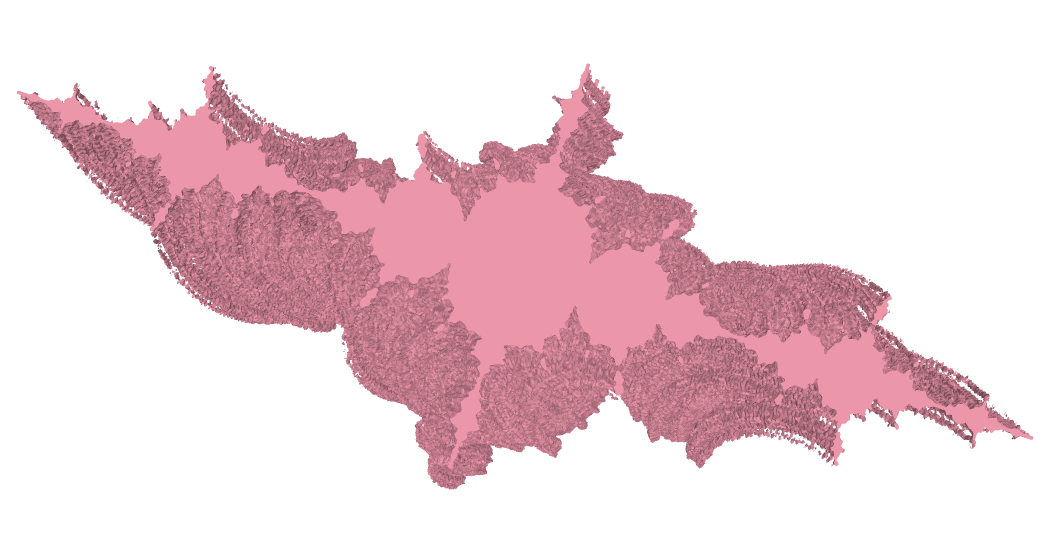}\hspace{10pt}
        \includegraphics[width=.3\paperwidth]{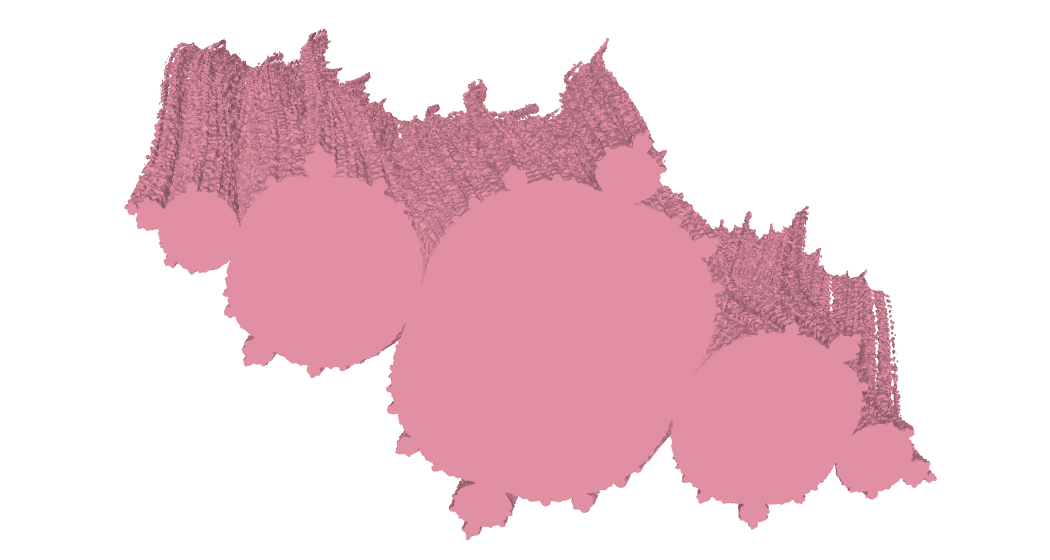}
    \end{center}
    The other hyperbolic components do not have boundaries that are so easy to parametrize.  Giarruso and Fisher explicitly parametrized the period 3 hyperbolic component \cite{period3}, and the parametrization is much more complicated than the two we have seen so far.  Beyond that, little has been worked out explicitly, leaving the question of how to correctly traverse the boundary as an interesting and difficult one.  In fact, the question of whether the boundary of the Mandelbrot set can be parametrized as a curve at all is related to the famous Mandlebrot Locally Connected (MLC) conjecture of Douady and Hubbard \cite{DHI},\cite{DHII}.  It probably isn't quite as difficult to find the boundary of a specific hyperbolic component, nevertheless, the problem of finding parametrizations of subsets of the boundary of the Mandlebrot set is connected to deep research questions in complex dynamics.\\
\subsubsection{Traversing Radii}
    The boundary is not the only interesting place to look.  Another direction (suggested by Caroline Davis), would be to consider the \textit{radii} of the hyperbolic components.  Indeed, each bulb of the Mandlebrot Set gives rise to Julia sets with different discrete invariants (such as the periods of the \textit{attracting cycles} \cite[Chapter 12]{Milnor}).  It is also known that each component is canonically holomorphic to a circle \cite{Milnor2}, so that each bulb has a canonically defined center and radii.  Tracing the radii between the center of one bulb to the center of another gives rise to a canonical deformation between two Julia sets with distinct discrete invariants, something that we are interested in 3D printing.  An important step is to obtain explicit parametrizations of these radii.  As with the boundaries, these are easy to find for the main cardioid, and the hyperbolic component of period 2 (which is already a disk), and can be extracted from \cite{period3} in the period 3 case.  Below is an image made by Walter Hannah \cite{WHannah} depicting these radii in the main cardioid.
    \begin{center}
        \includegraphics[width=.4\paperwidth]{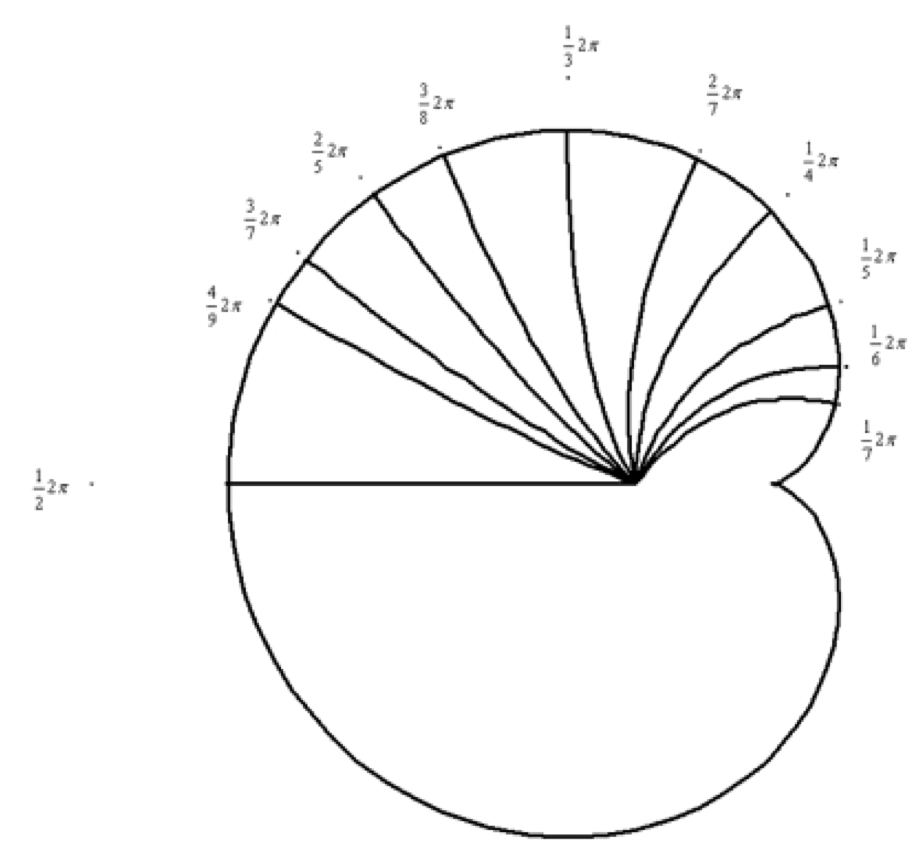}\cite{WHannah}
    \end{center}
\subsubsection{Related Objects in Holomorphic Dynamics}\label{CarolineUndergrad}
    We have only considered Julia sets associated to quadratic polynomials, but one can also iterate other functions$-$such as polynomials of higher degree or rational functions$-$and study their associated Julia sets.  An undergraduate research group at Indiana University, led by Caroline Davis, studied a family of polynomials whose associated Julia sets interpolate between the Apolonian and Sierpinski gaskets \cite{log(iu)}.  Below is a print of the associated static animation, printed by Bethany Mussman.
    \begin{center}
    \includegraphics[width=.6\paperwidth]{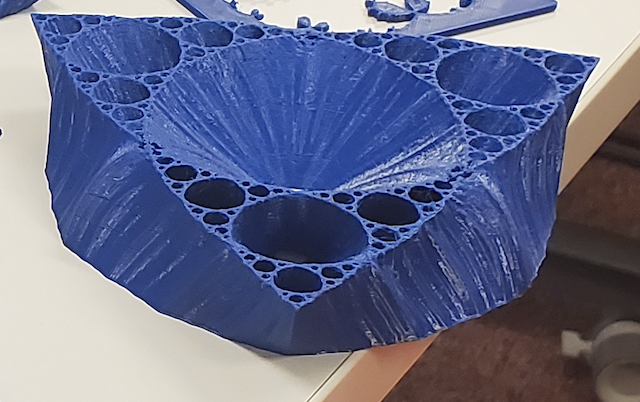}
    \end{center}
\subsection{Static vs. Computer Animations of Julia sets}
We are certainly not the first to visualize how Julia sets deform according to carefully chosen paths in the Mandlebrot set.  There are many computer animated visualizations of deforming Julia sets, including witnessing degenerations from connected Julia sets to Cantor sets, and exploratations near parabolic implosions.\footnote{In fact, the suggestion to 3D print parabolic implosions was made to the author by John H. Hubbard as he was demonstrating an animation he made of this phenomenon.}  The shift from computer to static animations brings with it advantages and disadvantages.  For example, computer animations give a more direct view of the internal structure, which can be hidden by the walls of the static animation, although this can be mostly remedied by adjusting the start and endpoints of the chosen path like we did in Section \ref{parabolicImplosions}.   Nevertheless, there is a completely novel perspective the static animation brings to the table, something that is missing from the computer animation: \textit{the side view}.  From the horizontal perspective we can observe how the different bulbs on the edge of these Julia sets move at different rates as the Julia set deforms, and these bulbs braid around each other in a way that seems to interpolate between the mathematical and the organic.  The following image gives a particularly good illustration of this phenomenon.
\begin{center}
    \includegraphics[width=.825\textwidth]{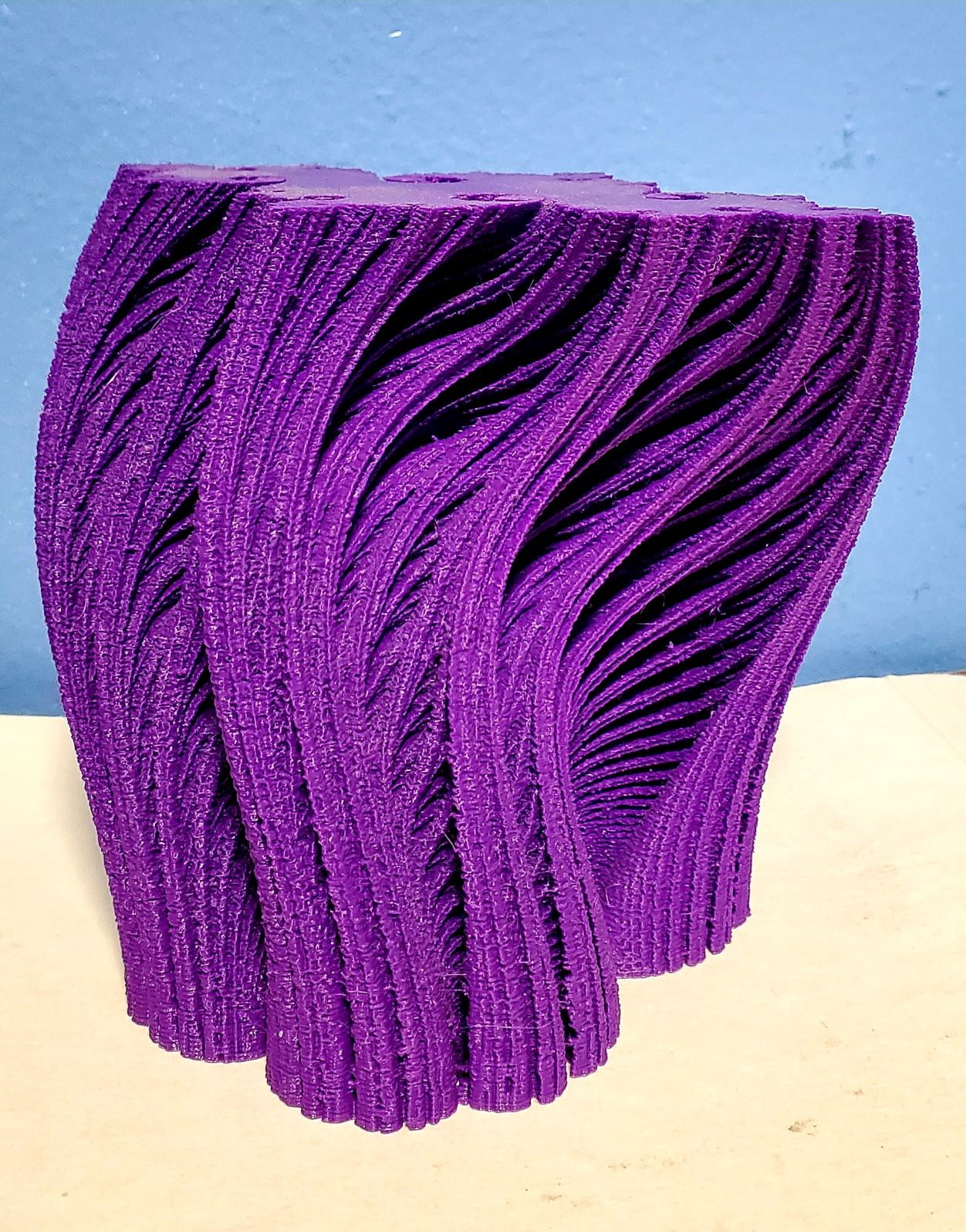}
\end{center}
\section{Creating a Static Animation}\label{How?}
When entering the world of 3D printing, it is easy to get overwhelmed by the large number of CAD programs and 3D modeling tools, many of which can be difficult to learn, especially as a beginner.  On the other hand, since static animations consist of stacks of 2-dimensional objects, one can hope to build a 3D static animation using purely 2-dimensional techniques.  In this section we describe a workflow that does exactly that,\footnote{This particular workflow was first shared with the author by Bernat Espigul\'e, and is how we created the static animations of filled Julia sets described in Section \ref{juliaSection}.} so that anyone who can make a 2-dimensional computer animation can begin designing 3D prints.  After describing the general process, we will then run through an explicit tutorial, building the static animation of deforming polar flowers from Section \ref{Theory}.  We intend this section to serve as an invitation for beginners to mathematical 3D printing, while also providing a new modeling technique to the experts.\\

The main idea is simple, and can be thought of as a sort of \textit{reverse tomography}.  While tomography is an imaging technique which extracts a collection of 2D cross sections from a 3D object, we begin with a collection of 2-dimensional frames, which we stitch together into a 3D model whose horizontal cross sections will be precisely the frames we begin with.\\

A summary of the workflow is the following: create a collection of 2-dimensional frames however you like, and upload them to a piece of software that can bond them together into a 3D model.  We use the open source software, UCSF Chimera \cite{chimera},\cite{chimeraSoftware}.  It has many uses, but we are most interested in its use as medical imaging software.  Roughly, MRI tomography takes a series of images of an organ at various depths inside the body, and these images can be fed to Chimera to reverse engineer a 3D model of the organ.  We can exploit this functionality to make our static animations.

\subsection{Overview}
Fix a collection of frames $B_t$, parametrized by values $t$ in some interval $I$.  We will describe a general process that can be followed to model the static animation
\[\cB = \{(b,t):b\in B_t\}\subseteq\bR^3,\]
as a 3D object file.  These steps assume the reader can already produce a 2-dimensional image of each frame $B_t$.
\subsubsection*{Step 1: Make a collection of frames}
Subdivide your interval into a reasonable number of frames.  Keep in mind that a common layer thickness for a 3D print is about .2mm, and the larger at home 3D printers allow a maximum height of no more than 300mm, so the resolution of your 3D printed object will probably not benefit from having more than 1500 frames (although industrial 3D printers can do better).  This gives a finite number of $t_i\in I$ which you will sample from.\\

For each $t_i$, produce a \textbf{PNG image} of $B_{t_i}$, coloring points $p\in B_{t_i}$ white, and points $q\notin B_{t_i}$ black.  This is the color scheme and file type that Chimera recognizes as cross sections of a 3D object.  Save the entire collection of images in a single folder and choose filenames so that the alphabetical ordering of the images corresponds to the correct ordering of the frames, this way they will be imported in the correct order during the next step.\\

You can create the frames of the animation any way you would like.  For the prints in Section \ref{juliaSection}, we fed a sampling of points on our parametrized path to Mathematica's JuliaSetPlot function \cite{Mathematica}, which we modified to color points in white and black as needed.\\

\textit{Warning: You may need to remove the alpha channel.}  The PNG graphics file format can store color images with an RGB color palette (storing numerical values for the red, green, and blue of each pixel), or an RGBA color palette, which has a fourth \textit{alpha} channel for transparency.  \textit{Chimera will reject PNG graphics with an alpha channel}, so  you may need to manually remove the alpha channel from your frames.  To do this we wrote a short script called \verb|alphaRemover.py|, included in the repository accompanying this paper \cite{github}.

\subsubsection*{Step 2: Convert the frames to a solid in Chimera}
From the Chimera home screen, choose \verb|File-->Open|.  Navigate to and highlight the entire set of images from Step 1, and select \verb|Open|.  A rough version of your model should become visible, as well as a \textit{volume viewer} menu with various settings.  There are 3 important settings to be aware of:
\begin{itemize}
    \item \textbf{Step:} A dropdown menu labelled \textit{step}, with choices of values 1,2,4,8, and 16 (4 is the default).  This roughly controls the vertical resolution.  \textit{Choose 1} (the highest vertical resolution).  Your model should immediately look much sharper.
    \item \textbf{Level:} A slider to select the \textit{level}, which is a number between 0 and 255.  This roughly controls the tolerance for edge detection.  We often have to play around with this a bit.
    \item \textbf{Style:} Radio buttons to select either surface, mesh, or solid.  We always choose \textit{surface}.
\end{itemize}
Once you have found settings to your liking, you can export the model.  Select \verb|File-->Export Scene...|  You can choose from a number of file types.  The standards for 3D printing are OBJ or STL.  Depending on the size of the image and number of frames, \textit{this can take a long time}.  For the particularly detailed animations of Section \ref{juliaSection} we needed to let it run for several hours.\footnote{For reference, the computation was done on a modestly powered Dell XPS 15 7590 with an Intel Core i7-9750H (12MB Cache, up to 4.5 GHz, 6 cores).}

\subsubsection*{Step 3: Slice and Print}
Your model should now be ready to open up in a slicer (for example \textit{Cura} \cite{Cura}), and then sliced and sent to a 3D printer.  The full process depends on your 3D printer, and is beyond the scope of this article, but we'd like to include a few tips to help avoid failure at this step.  We point to a few resources for slicing and printing 3D models at the end of this subsection.\\

\textit{Avoid steep slopes and overhangs.}  An extrusion based 3D printer will print your static animation by printing each frame on top of the last, using the previously printed frames as a support base for the next ones.  If the frames $B_t$ vary too wildly with $t$, the printer may end up trying to extrude part of a frame with no supporting structure underneath.  These unsupported layers of plastic building up will lead to unsatisfactory results and failed prints.  As a rule of thumb, most 3D printers can cleanly print overhanging structures with an angle up to somewhere between 45 and 60 degrees, so if you would like to print a wildly deforming family, you may need to do some experimentation to get something that works.\footnote{Alternatively, you can employ the use of \textit{support material}, extra printed material which holds up excessive overhangs and can be removed after.}\\

\textit{Lay your model flat.}  Chimera won't automatically export a model which lies flat, instead exporting the model aligned to the perspective in the viewing window.  Sending the tilted model to an extrusion based 3D printer will likely cause issues with adhesion and support, and could result in a failed print.  This can be avoided if you don't adjust the view of the model in Chimera's viewscreen between importing your frames and exporting the STL.  Otherwise you must rotate your exported model so that it lies flat on the buildplate.  This can be done in Cura using the \textit{Lay Flat} command.\\

\textit{More resources for the printing process.}  Although beyond the scope of this paper, there are many places to turn for a thorough introduction to the use of slicers and 3D printers to turn a digital model like your static animation into a physical model.  For example, the author's lecture during the 2022 Joint Mathematics Meetings mini-course \textit{3D Printing: Challenges and Applications} was a tutorial on slicing and printing, and the slides are publicly available \cite{mySlides}.  There are also many detailed written tutorials (for example \cite{printerlyTutorial}) and youtube tutorials (for example \cite{youtubeTutorial}).

\subsection{Tutorial}\label{tutorial}
We now demonstrate this process step-by-step with an example, constructing the static animation of deforming polar flowers introduced in Section \ref{Theory}.  The reader is invited to follow along on their own machine.  To this end we have provided code snippets within this document, as well as the complete collection of code and files in a public Github repository \cite{github}.  For more detail on where to find what you need to follow along, see Section \ref{where} below.\\

\begin{remark}
In Section \ref{preciseTheory} we computed a complete parametrization of this particular deformation space in $\bR^3$.  With this parametrization there are many more direct ways to produce a 3D model (for example, plotting in Mathematica or Sagemath and then exporting as an STL).  Nevertheless, one often encounters deformation spaces which do not admit such a parametrization (for example, those in Sections \ref{treeSection} and \ref{juliaSection}).  One should think of this tutorial as an introduction to this particular process, rather than the most efficient way to make this particular model.
\end{remark}

The frames of the static animation in this tutorial are polar curves whose equations we computed in Section \ref{preciseTheory} to be $r = 2 + t\cos(5\theta + 2\pi t)$ for values $t$ in the interval $I = [0.2,1]$.  Since we are making a 3-dimensional solid, we will replace these curves with the regions they bound, so that our resulting object will be filled rather than hollow.  The frames we will use will therefore be the following sets
\[B_t = \{(r,\theta): r\le 2+t\cos(5\theta + 2\pi t)\}\subseteq\bR^2.\]
For example, $B_1$ is now the entire shaded region below, rather than just the boundary curve.
\begin{center}
    \[\includegraphics[width=.25\textwidth]{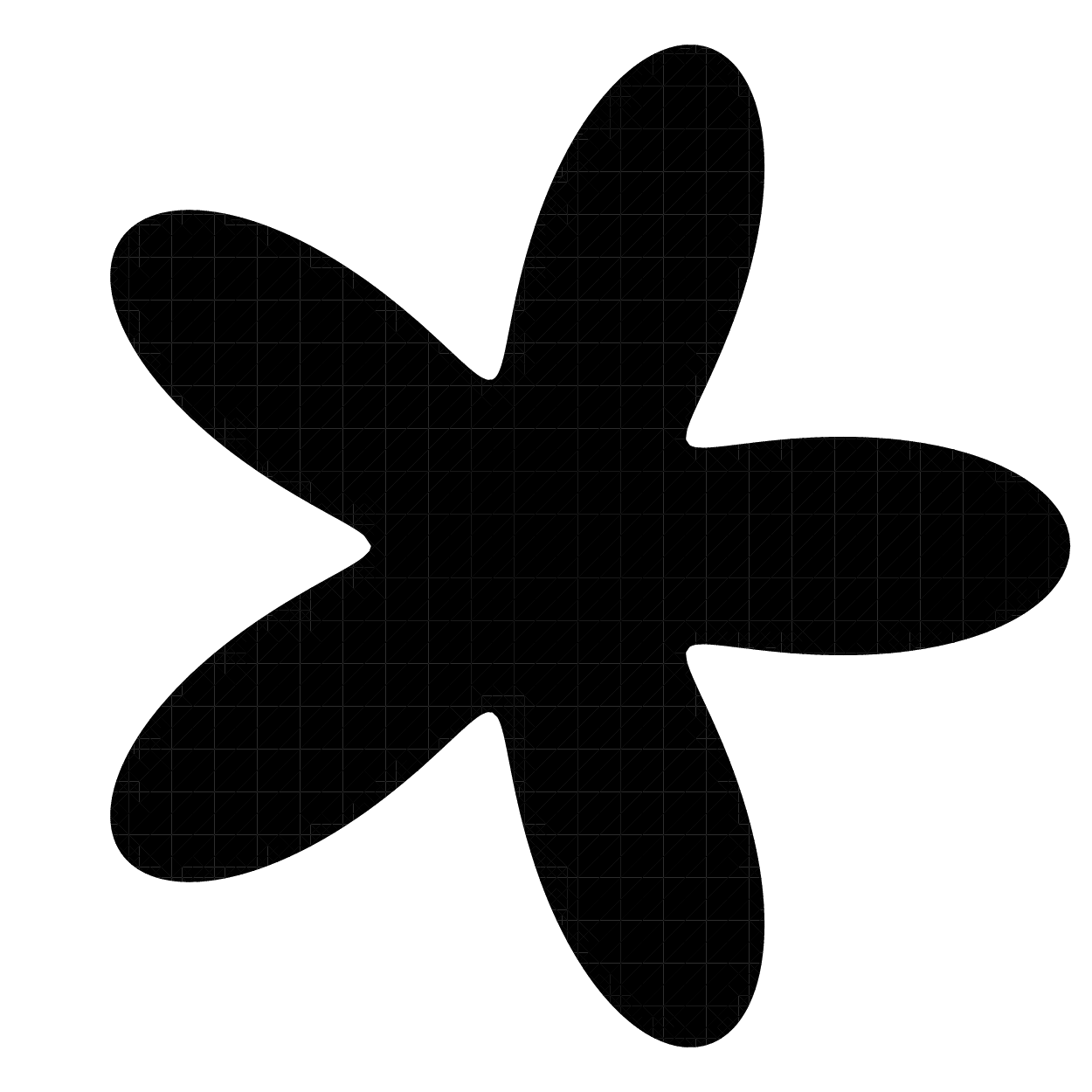}\]
\end{center}
We now run through all the steps to model the static animation $\cB = \{(b,t):t\in I,b\in B_t\}$ as a 3D object file.
\subsubsection*{(Tutorial) Step 1}
We subdivide our interval to create 161 frames, so that $t$ progresses between $0.2$ and $1$ at steps of $.005$ units.  Our list of $t$ values is therefore
\[T = \{.2, .205, .21, .215,\cdots,.990,.995,1\}.\]
We need to make a PNG image for each frame, drawing each region $B_t$ in white, against a black background.  Let's walk through a way to do this in python.  We need the \verb|matplotlib| and \verb|numpy| libraries.
\begin{python}
  import matplotlib.pyplot as plt
  import numpy as np
\end{python}
Chimera expects white images against a black background, so we first set the color of the background we are drawing against.
\begin{python}
  plt.style.use('dark_background')
\end{python}
For each frame, we'd like to draw the polar curve $r = 2t\cos(\theta+2\pi t)$.  To do this, \verb|mathplotlib| will sample a number of points and interpolate.  We use numpy's \verb|linspace| command to create the list of $\theta$ values to sample from.
\begin{python}
  NUM_PTS = 1000
  theta = np.linspace(0,2*np.pi,NUM_POINTS)
\end{python}
We then calculate the corresponding $r$ values to finish computing the list of points.  Notice we set $t=1$ in this snippet, but this will change from frame to frame.
\begin{python}
  t = 1 #This varies on each frame
  r = 2 + t*np.cos(5*theta + 2*np.pi*t)
\end{python}
Now whats left is to let \verb|mathplotlib| draw and fill the curve.
\begin{python}
  plt.polar(theta,r,alpha=0) #Draws the curve
  plt.fill_between(theta,0,r,color='white') #Fills the curve
  plt.axis('off') #We just want the flower, no axes
\end{python}
Here is the output with $t$ values 1 and 1/2 respectively.\\

\begin{center}
    \includegraphics[width=.3\paperwidth]{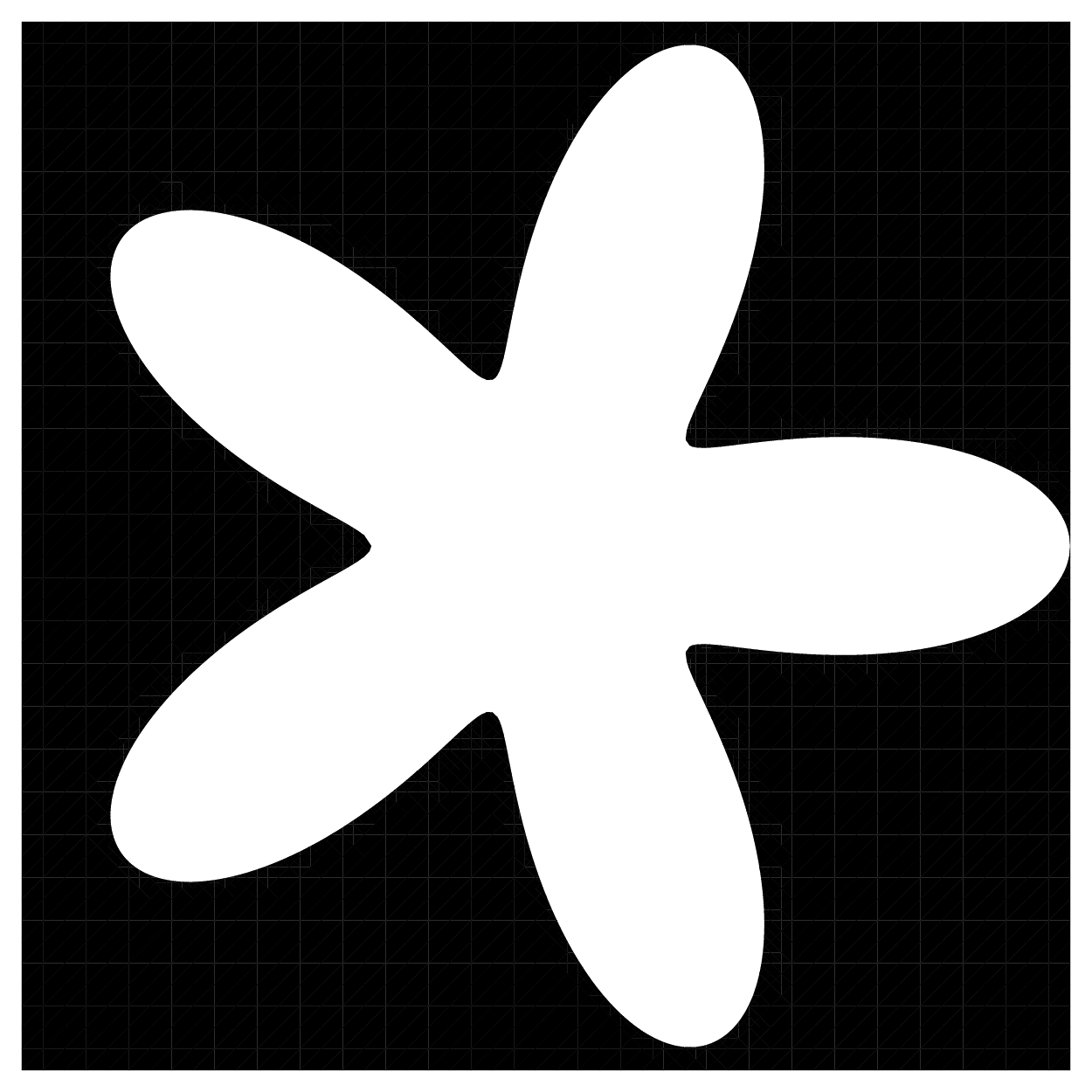}\hspace{10pt}\includegraphics[width=.3\paperwidth]{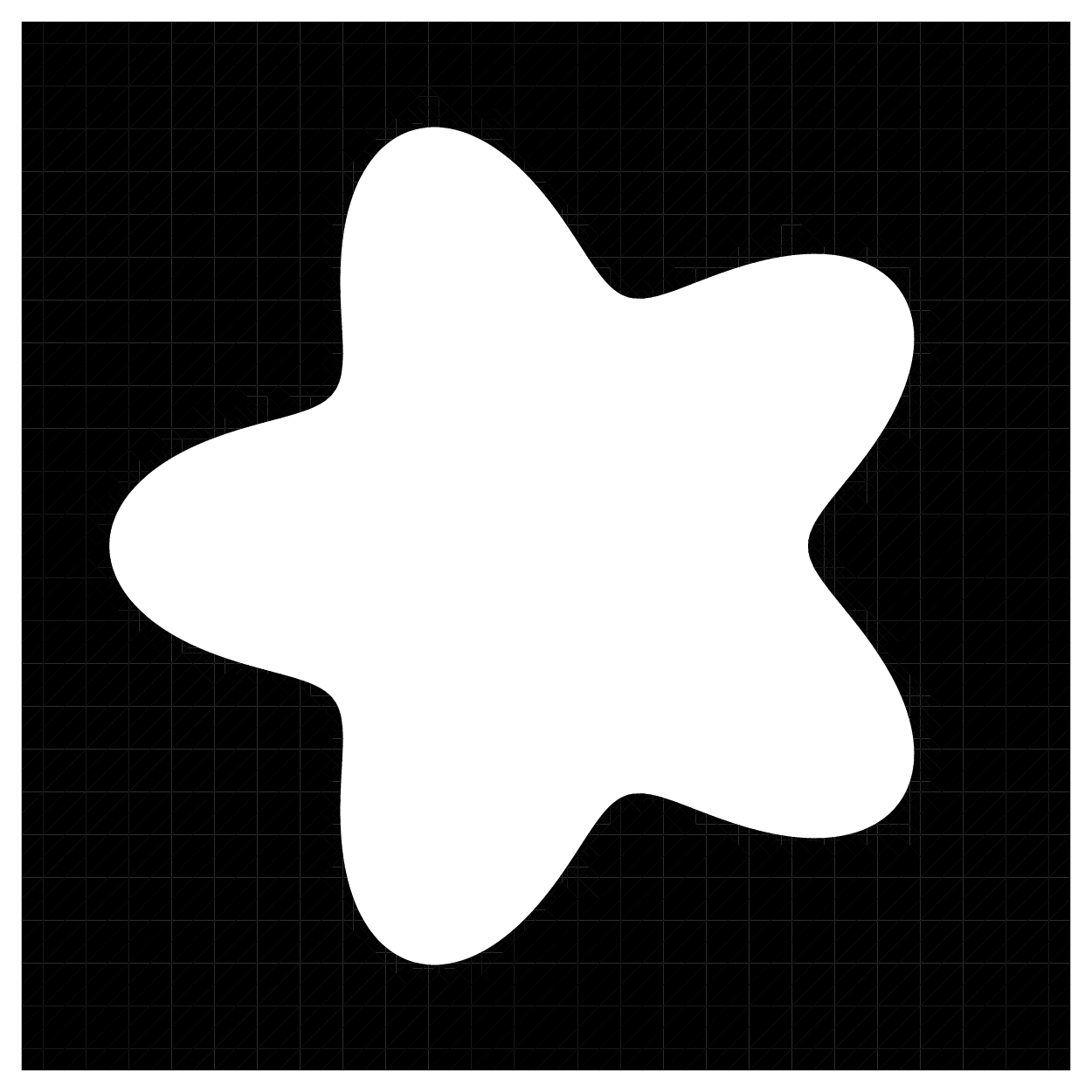}
\end{center}

To create the entire stack of frames we loop through each value of $t\in T$ and export a PNG built using the the code above, choosing filenames in alphabetical order so that they are arranged correctly when imported by Chimera.  We can do this by looping through the last two snippets of the code above and saving the output for each frame.
\begin{python}
NUM_FRAMES = 161 #Set the number of frames

for i in range(NUM_FRAMES):
  t = 0.2 + .005*i
  r = 2 + t*np.cos(5*theta + 2*np.pi*t) #The boundary curve
  plt.polar(theta,r,alpha=0) #Draws the curve
  plt.fill_between(theta,0,r,color='white') #Fills the curve
  plt.axis('off') #We just want the flower, no axes

  #Now export and clear the plot.
  imageName = 'flowerStack/flowers_'+str(1000+i)+'.png'
  plt.savefig(imageName) #Export the image
  plt.cla() #clear the plot before looping back
\end{python}
The images are now saved as \verb|flowers_1000.png| through \verb|flowers_1160.png|.  The reason we started counting at 1000 is because Chimera reads alphabetically in such a way that 10 comes before 9 because it starts with a 1.  In particular, if we started counting at 0 Chimera would arrange the imported frames in an incorrect order during Step 2.  We now have our frames images all saved as PNG images.
\begin{center}
    \includegraphics[width=.45\paperwidth]{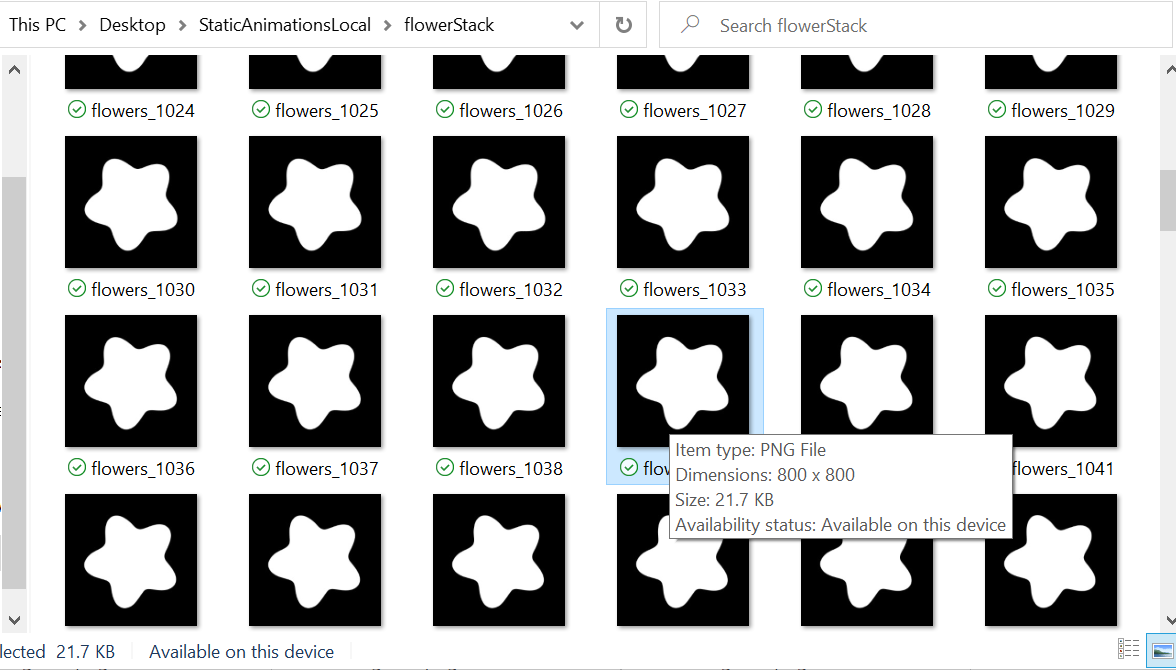}
\end{center}
In our accompanying github repository \cite{github} we include a python program \verb|flowerPlot.py| for creating the entire stack of frames, as well as an implementation in Mathematica, \verb|flowerPlot.nb|.
\subsubsection*{(Tutorial) Step 2}
Open Chimera, select browse, and highlight all of the PNG images from step 1.
\begin{center}
    \includegraphics[width = .35\paperwidth]{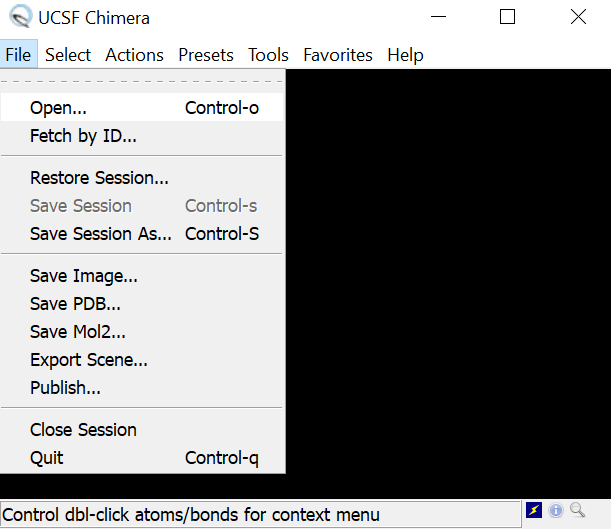}\hspace{20pt}\includegraphics[width=.25\paperwidth]{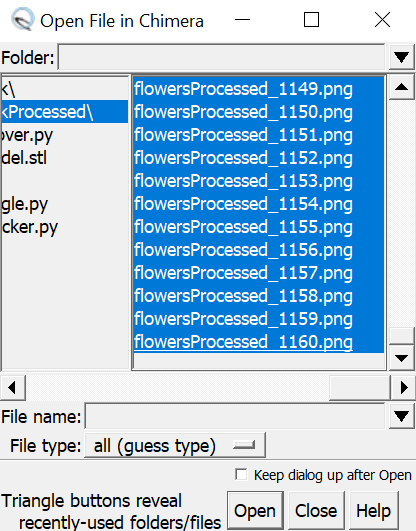}
\end{center}
Click open.  Chimera stitches all of the polar flowers into a single solid.
\begin{center}
    \includegraphics[width = .6\paperwidth]{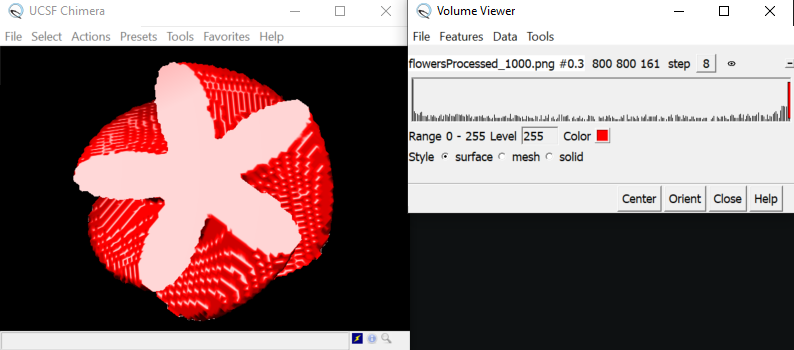}
\end{center}
It looks a bit rough.  By changing \textit{step} to 1 and adjusting the level we arrive at a much finer model.  We also choose \textit{surface} in the style setting (if it isn't already selected).
\begin{center}
    \includegraphics[width = .6\paperwidth]{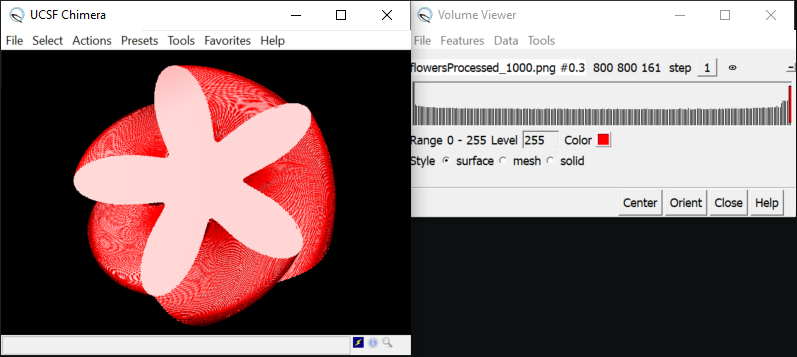}
\end{center}
This looks much better.  Here are a few other points of view of the model Chimera produced.
\begin{center}
    \includegraphics[width=.3\paperwidth]{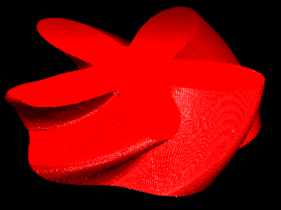}\hspace{20pt}\includegraphics[width=.3\paperwidth]{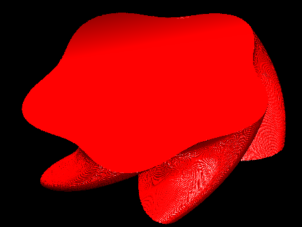}
\end{center}
We now export the model from Step 2 as an STL.  To do this we merely export the scene, choosing the STL file extension.  We should make sure to export the model before rotating it at all in the viewscreen.  This way it will be oriented correctly when opened in a slicer, lying flat against the buildplate (see the Step 3 outline above).
\begin{center}
    \includegraphics[width = .45\textwidth]{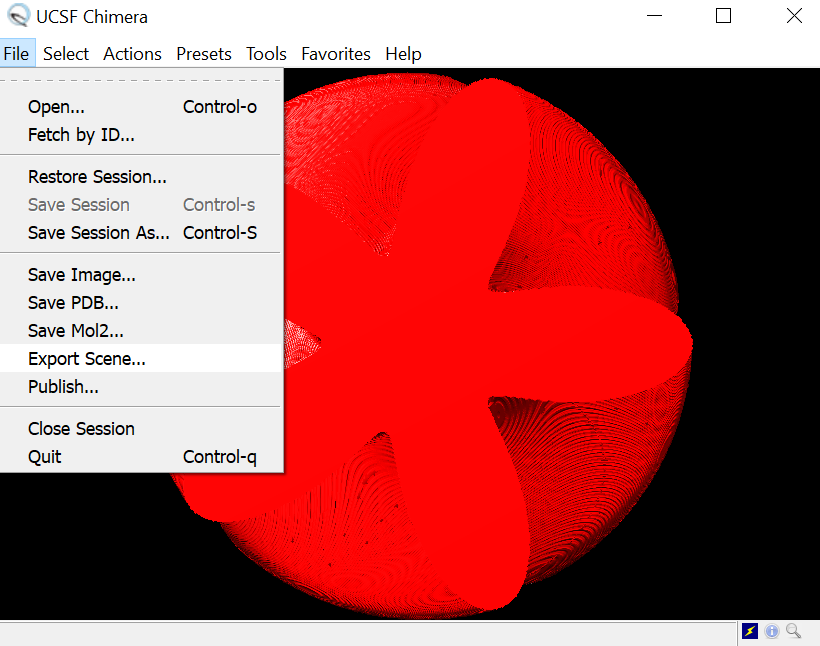}
    \hspace{10pt}
    \includegraphics[width = .45\textwidth]{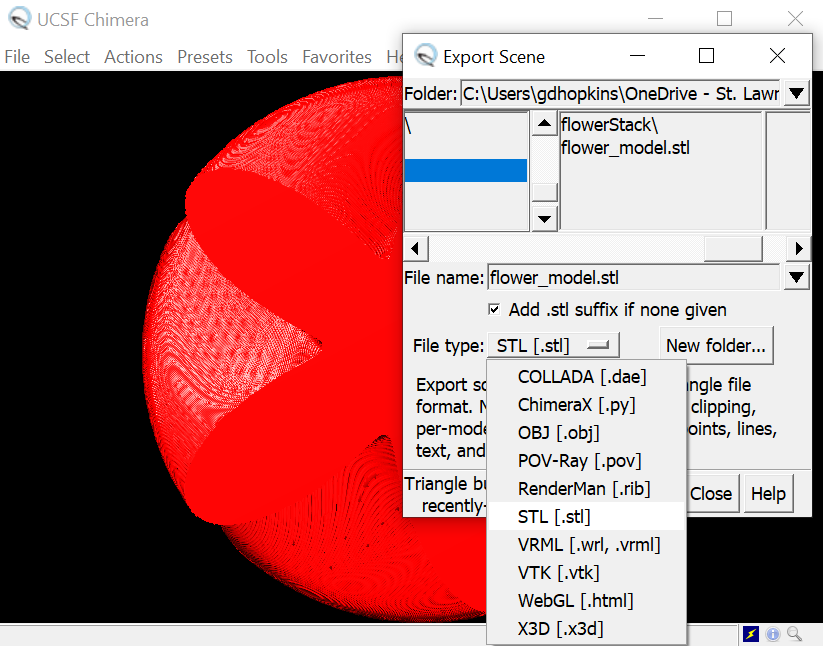}
\end{center}
For more complicated animations this can take a long time$-$for example, if you need more frames, or a higher resolution to capture the complex behavior you are interested in.  This animation, though, is relatively simple, and on the author's current machine\footnote{HP ProBook 445 G8, AMD Ryzen 5 5600U with Radeon Graphics 2.30 GHz, RAM: 16 GB} it took about 3 minutes to export.\\

\textit{Warning.} When importing your frames to Chimera at the beginning of this step, you may have encountered the following error message:
\begin{center}
  \includegraphics[width=.55\textwidth]{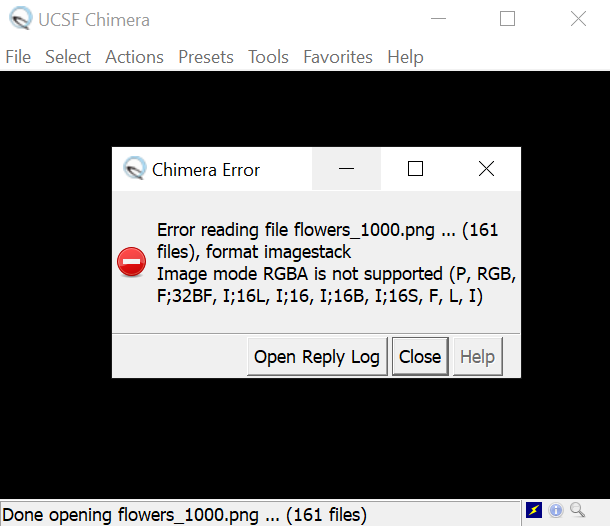}
\end{center}
This is because the frames were saved as PNG files \textit{with an alpha layer}.  To fix this, you must remove the alpha layer.  One can do this in Python using the Python Imaging Library (PIL).
\begin{python}
from PIL import Image
image = Image.open(fileName).convert('RGB') #conversion step
image.save(fileName)
\end{python}
In the snippet above, \verb|fileName| should be a path to your image.  In the github repository accompanying this manuscript \cite{github} we include a script, \verb|alphaRemover.py|, which can loop through your entire image stack and remove the alpha layer from each frame.  The program \verb|flowerPlot.py| accompanying this tutorial already removes the alpha layers from your stack, so this step will only be necessary if you created your image stack elsewhere.
\subsubsection*{(Tutorial) Step 3}
The model is now ready to be opened in slicing software (we use Cura \cite{Cura}), sliced, and then sent to a printer.
\begin{center}
    \includegraphics[width = .65\paperwidth]{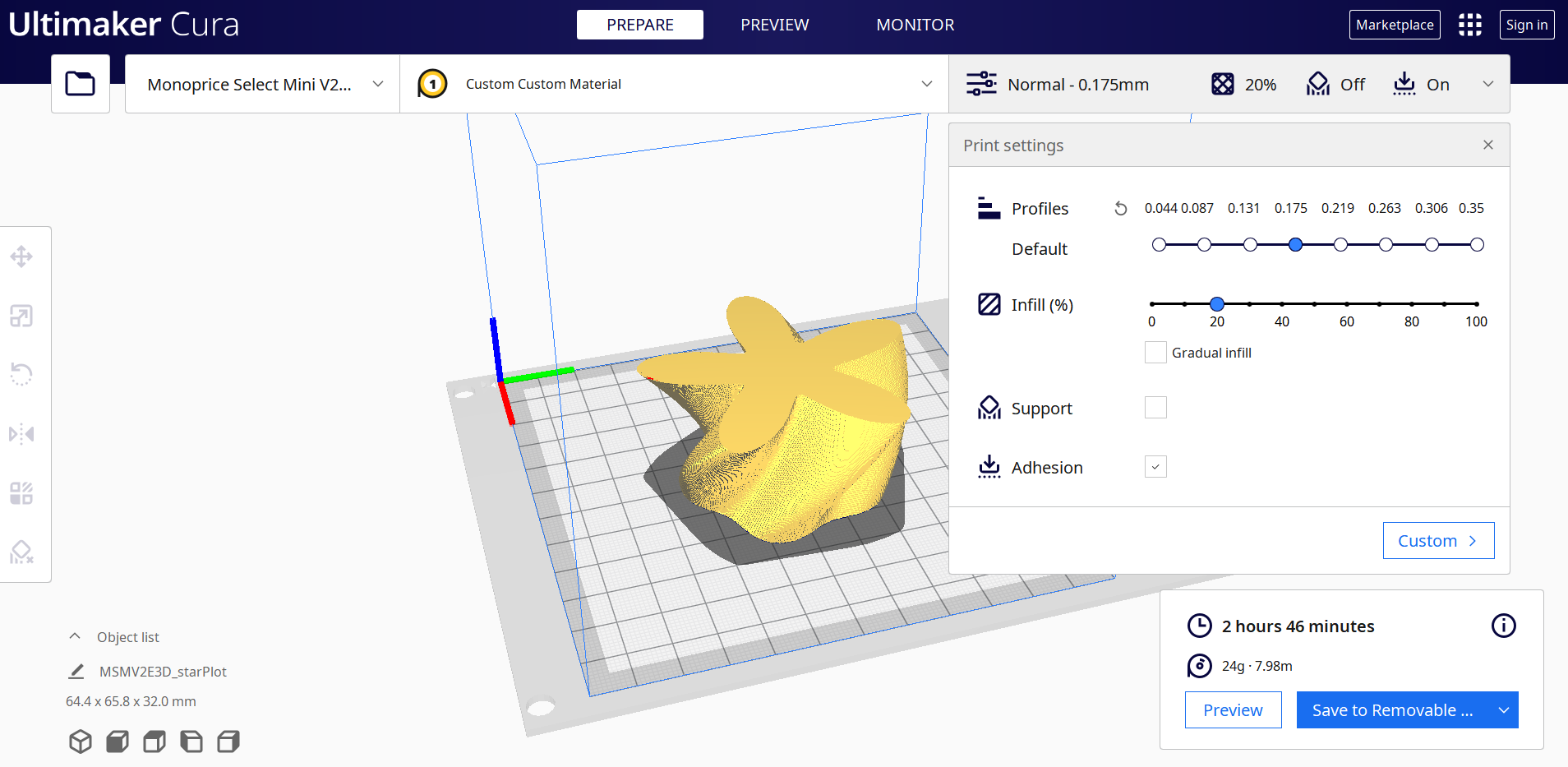}
\end{center}
\begin{center}
    \includegraphics[width=.65\paperwidth]{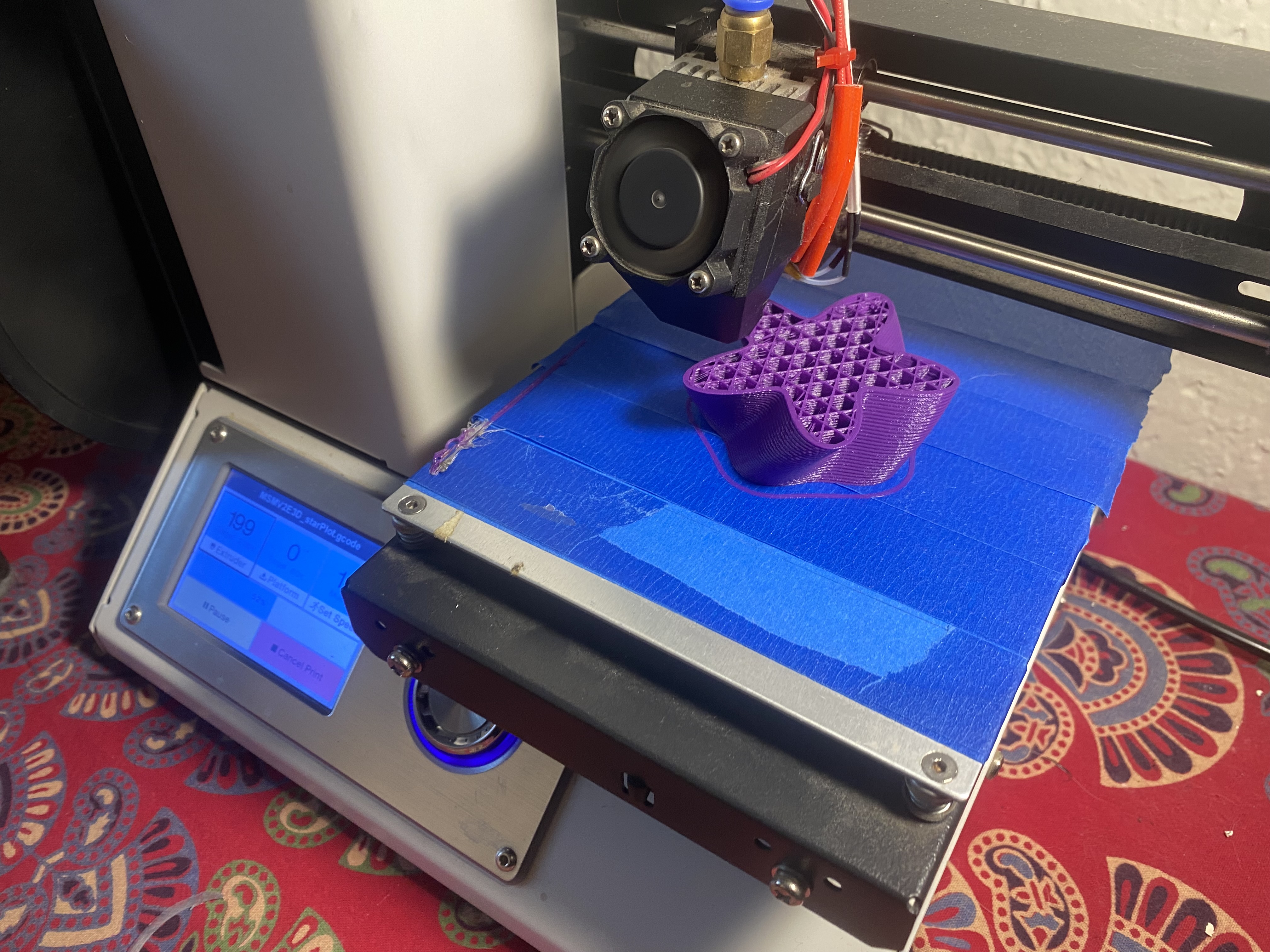}
\end{center}
\begin{center}
    \includegraphics[width=.3\paperwidth]{images/flowerSide.jpg}\hspace{10pt}
    \includegraphics[width=.3\paperwidth]{images/flowerTop.jpg}
\end{center}
\subsection{Where do I get...}\label{where}
Our intent with this tutorial is to have the reader follow along on their own machine.  This way they learn the workflow and can integrate it into their own mathematical practice.  In order to follow along, the reader will need the following tools.
\begin{itemize}
    \item \textbf{Chimera:} At the time of writing of this document, UCSF Chimera \cite{chimeraSoftware} is available for free download here: \url{https://www.cgl.ucsf.edu/chimera/}.
    \item \textbf{Cura:} At the time of writing of this document, Untimaker Cura \cite{Cura} is available for free download here: \url{https://ultimaker.com/software/ultimaker-cura}.
    \item \textbf{Tutorial Files:} All the files generated while making this tutorial are available at our public github repository \cite{github}.  This includes:
    \begin{itemize}
        \item A python program generating the image stack.
        \item A Mathematica notebook generating the image stack.
        \item A python program for removing the alpha layer from a PNG stack.
        \item The PNG images that serve as the frames of the animation.
        \item The STL exported from Chimera.
    \end{itemize}
 \end{itemize}
    If the reader would only like to follow along for one step, for example, stitching together the image stack in Chimera, they are welcome to download the frames from \cite{github} and experiment in Chimera themselves.
\section{Future Directions}\label{conclusion}
This article only begins to scratch the surface of what the author believes to be a vast array of possibilities for 3D printed static animations in mathematical illustration and research. Indeed, with how ubiquitous the notion of deformations are in mathematics, and how useful animations on a screen have proved to be, we look forward to exciting new directions for this technique of illustration and experimentation.  We hope that with the techniques developed in Section \ref{How?}, the reader will have the tools necessary both to incorporate 3D printing into their practice, and to find new applications for static animations.  We'd like to conclude with a couple of potential next steps.

\subsection{Exploring the Mandlebrot set}
As we saw in Section \ref{juliaMore}, the field of holomorphic dynamics offers a rich menagerie of static animations to investigate, including many stemming from paths within the Mandlebrot set. The process for creating static animations from paths in the Mandlebrot set is already established (and can be accessed in our public github repository \cite{github}), so all one needs are parametrized paths and to begin to explore.
\subsection{Chaos and the Double Pendulum}
The double pendulum \cite{chaos} is an example of a chaotic system.  In this system, attached to the mass at the end of a pendulum there is a second pendulum with another mass.
\begin{center}
    \includegraphics[width=.15\paperwidth]{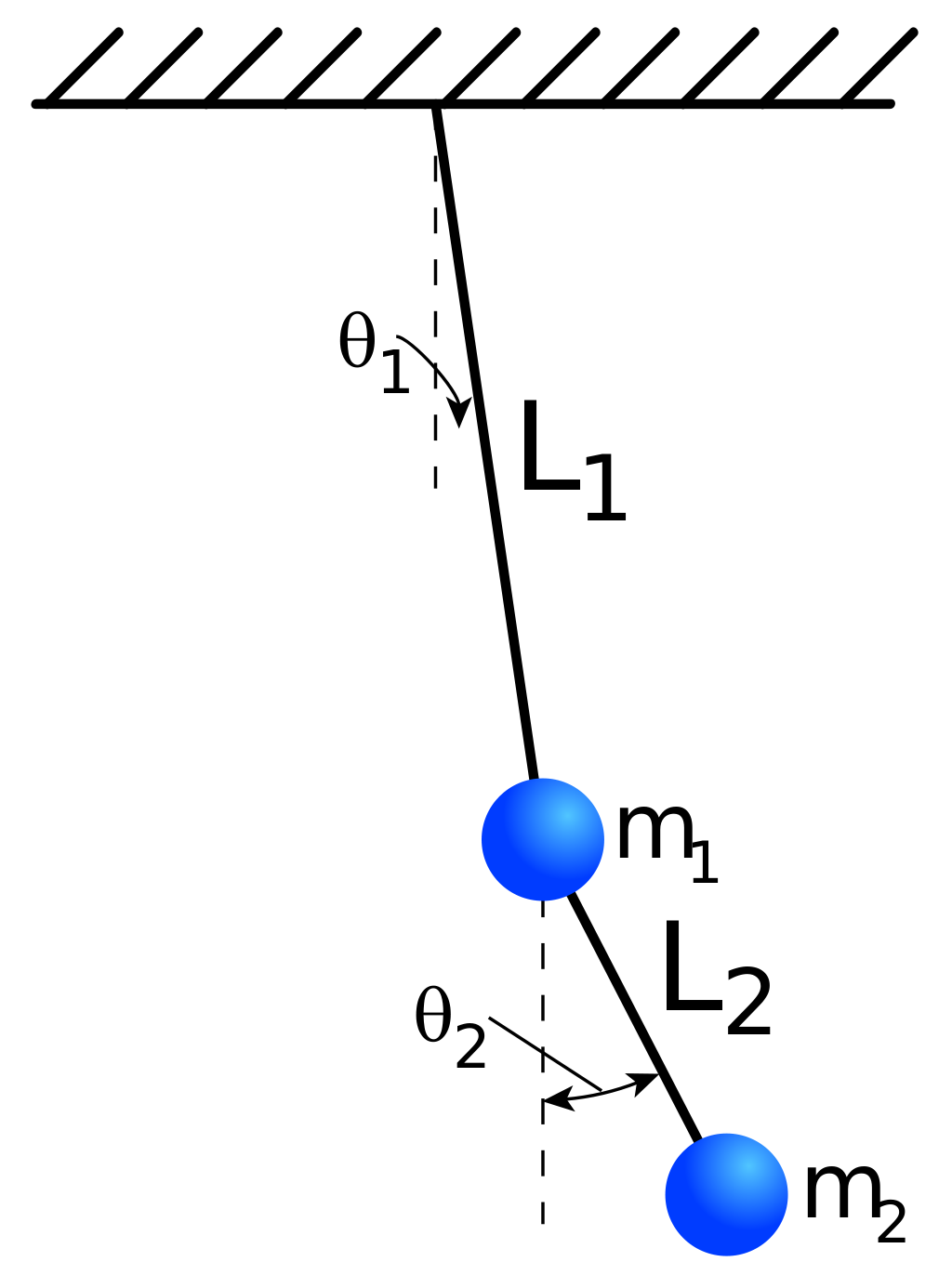}
    \cite{wikipedia}
\end{center}
This is a dynamical system governed by a set of ordinary differential equations which is \textit{chaotic}, that is, it is highly sensitive to its initial conditions.  Animations and experiments with double pendulums show that, even with two very similar starting points, the behavior of the moving pendulums will quickly diverge quite radically.  We imagine turning the computer animation into a static animation, building the subsequent positions of the double pendulum on top of the last.  We have begun such an exploration with the help of Daniel Piker using his physics engine, Kangaroo \cite{kangaroo}.
\subsection{And more...}
We invite the reader to include static animations in their own investigations.   Any place that a 2D animation has been useful, explore adapting it to a 3D static animation.  We hope that, with the gentle learning curve and connections to many mathematical disciplines, folks will find this to be an accessible avenue to the world of 3D printing in mathematics.  Indeed, as we showcased in Section \ref{CarolineUndergrad}, we have begun to see researchers lead undergraduate research investigations which use these techniques \cite{log(iu)}.\\

There are many fascinating applications of 3D printing to exploring and illustrating mathematics, and we hope that static animations will both become a part of that canon, and also serve as an invitation to the broader world of 3D printing in math.

\bibliography{bib}{}
\bibliographystyle{acm}

\appendix
\section*{Appendix: Deformations in Algebraic Geometry}
The author's initial motivation to create 3D printed static animations of continuously deforming mathematical objects came from encountering these types of constructions in algebraic geometry, and more specifically within deformation theory.  Here, one often will think of parameters of variation as a spacial variables and produce a deformation space of higher dimension to which witnesses the interpolation.  Indeed, if $C$ and $C'$ are varieties which can deform to one another, exhibiting them both as living within some larger deformation space $D$ gives a geometric framework wherein information can be passed between them via the geometry of $D$.  As this technique was the authors inspiration to begin creating static animations, we'd like to use this appendix to describe this connection in some more detail, using examples of flat families of algebraic curves.\\

Suppose you are studying a family of curves $C_t$ parametrized by $t\in k$ where $k$ is a field (little is lost imagining $k$ to be $\bC$ or $\bR$).  Suppose further that each curve is the zero set in $k^2$ of a polynomial equation in 2 variables $x$ and $y$.  The key observation is that this $t$ can be introduced as a new \textit{spacial} variable, defining an equation, now in $x,y,$ and $t$, whose zero set in $k^3$ we will call $D$.  Then $D$ comes equipped with a projection $\pi:D\to k$, which plucks out the $t$-coordinate, and each member of the original family of curves can be recovered as the fibers of this projection:
\[C_t = \pi^{-1}(t)\subseteq D.\]
The remarkable thing about this viewpoint is that relationships between the various curves $C_t$ in the family can be understood in terms of the global geometry of the deformation space $D$, which can be studied as an algebraic variety in its own right. For example, one can use the geometry of $D$ to show that the genus of the curves in such a family is constant \cite[Corollary III.9.13]{Hartshorne}.\footnote{The genus of a curve can be defined over any field,  but if $k=\bC$, then the curve lives as a subspace of 1 complex dimenion in $\bC^2$ and can therefore be thought of as a Riemann surface, allowing one to determine genus in the usual way} To connect with our 3D printing philosophy, the 3D printed static animation plays the role of the total deformation space $D$, which can be viewed as a vertical stack of the $C_t$ at height $t$.  Much like how the global geometry of $D$ provides new insights into the connections between the parametrized curves $C_t$, this paper begins to ask how the global structure of a 3D printed static animation reflects the connections between the parametrized objects it contains.\\

We give one explicit example of the deformation space of a family of curves.  To do this let $k=\bR$ and begin by considering the hyperbola $xy = 1$.
\begin{center}
\begin{tikzpicture}
  \draw[->] (-3, 0) -- (3, 0) node[right] {$x$};
  \draw[->] (0, -3) -- (0, 3) node[above] {$y$};
  \draw[domain=.33:3, smooth, variable=\x, blue] plot ({\x}, {1/\x});
  \draw[domain=-3:-.33, smooth, variable=\x, blue] plot ({\x}, {1/\x});
\end{tikzpicture}
\end{center}
Now let us perturb the equation, considering $xy=\frac{1}{2}$, or $xy=\frac{1}{4}$.
\begin{center}
\begin{tikzpicture}
  \draw[->] (-3, 0) -- (3, 0) node[right] {$x$};
  \draw[->] (0, -3) -- (0, 3) node[above] {$y$};
  \draw[domain=.17:3, smooth, variable=\x, blue] plot ({\x}, {.5/\x});
  \draw[domain=-3:-.17, smooth, variable=\x, blue] plot ({\x}, {.5/\x});
  \draw[domain=.08:3, smooth, dashed, variable=\x, red] plot ({\x}, {.25/\x});
  \draw[domain=-3:-.08, smooth, dashed, variable=\x, red] plot ({\x}, {.25/\x});
\end{tikzpicture}
\end{center}
We can continue decreasing this constant term, considering the curve $C_t$ given by the equation $xy=t$ for $t=\frac{1}{8},\frac{1}{16},...$.
\begin{center}
\begin{tikzpicture}[samples = 50]
  \draw[->] (-3, 0) -- (3, 0) node[right] {$x$};
  \draw[->] (0, -3) -- (0, 3) node[above] {$y$};
  \draw[domain=.04:3, smooth, variable=\x, blue] plot ({\x}, {.125/\x});
  \draw[domain=-3:-.04, smooth, variable=\x, blue] plot ({\x}, {.125/\x});
  \draw[domain=.02:3, smooth, dashed, variable=\x, red] plot ({\x}, {.0625/\x});
  \draw[domain=-3:-.02, smooth, dashed, variable=\x, red] plot ({\x}, {.0625/\x});
\end{tikzpicture}
\end{center}
As $t$ approaches zero the curves appear to be converging to the axes.  And indeed, $xy=0$ if and only if $x=0$ or $y=0$, so that $C_0$ is precisely the union of the $x$ and $y$ axes.  We have therefore constructed a family of curves $C_t$, parametrized by $t\in\bR$, which continuously interpolate between hyperbolas and a union of 2 lines.  Of course, nothing is stopping us from treating $t$ as a \textit{spacial} variable.  That is, we can consider the the surface in $\bR^3$ defined by the equation $xy=z$.
\begin{center}
    \includegraphics[width=.3\paperwidth]{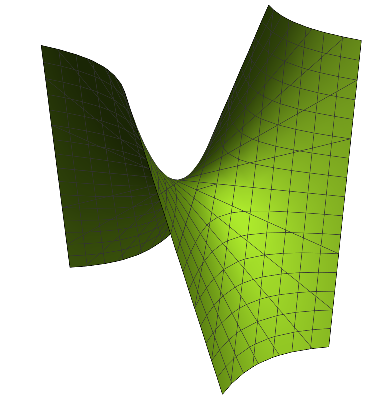}
\end{center}
The horizontal cross section at height $z = t$ is precisely the curve $C_t$ given by the equation $xy=t$, and if one considers the associated contour lines, one recovers the picture of the deformation of hyperbolas to the planar axes.
\begin{center}
    \includegraphics[width=.4\paperwidth]{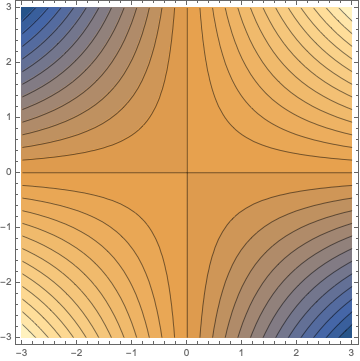}
\end{center}
In this way, the surface $xy=z$ is a higher dimensional space which witnesses the continuous deformation of hyperbolas to a union of 2 lines.\\
\end{document}